\documentclass[a4paper,12pt,oneside]{article}

\usepackage{graphicx}

\usepackage{latexsym}
\usepackage[english]{babel}
\newtheorem{theorem}{Theorem}
\newtheorem{lemma}[theorem]{Lemma}
\newtheorem{proposition}[theorem]{Proposition}

\newtheorem{definition}[theorem]{Definition}
\newtheorem{corollary}[theorem]{Corollary}
\newtheorem{remark}[theorem]{Remark}
\newtheorem{axiom}[theorem]{Axiom}

\begin{document}
\title{Orbit complexity, initial data sensitivity and weakly chaotic dynamical systems }

\author{Stefano Galatolo
\\
Dipartimento di Matematica  \\
and \\
 Centro Interdisciplinare per lo Studio dei
Sistemi Complessi,\\
 Universit\`a di Pisa \\
 Via Buonarroti 2/a, Pisa, Italy \\
e-mail: galatolo@dm.unipi.it}

\maketitle
\abstract {We give a definition of generalized indicators 
of sensitivity to initial conditions and orbit complexity (a measure of the 
information that is necessary to describe the orbit of a given point).
The well known Ruelle-Pesin  and Brin-Katok theorems, combined with
Brudno's theorem give a relation between initial data sensitivity
and orbit complexity that is generalized in the present work.
The generalized relation implies that the set of points
where the sensitivity to initial conditions  is more than exponential in all directions is a 0 dimensional 
set.
The generalized relation
is then applied to the study of 
an important example of weakly chaotic dynamics: the  Manneville map.}

\def\stackunder#1#2{\mathrel{\mathop{#2}\limits_{#1}}}

\section{Introduction}

%!@***da mettere dentro forse :
%Si fanno le cose senza misura invariante anche perche ci interessano i sistemi
%non all equilibrio, che vengono fuori dall autoorganizzazione e danno 
%lougo a processi non strazionari.@

When we think about a chaotic system there are two things that we expect to happen: 

\begin{itemize}
\item the behavior of the system is unpredictable and complex to be described

\item  small differences  on initial conditions leads to big differences in the evolution of the system.
\end{itemize}

A rigorous measure of initial conditions sensitivity is not difficult  to formalize, it leads to the definition of Lyapunov exponents (see e.g. \cite{pesin}) or to the more general Brin-Katok \cite{brinkatok} local entropy.

A  measure of the complexity of the behavior of a system is less simple 
to  formalize and  in the case of dynamical systems it was given by Brudno \cite{brud}. A set of strings is associated by a certain construction to the orbit of a point   and then the complexity of the orbit is defined by the Algorithmic Information Content  (section 2 ) of the associated strings.
The complexity of an orbit is then a measure of the amount of information 
that is necessary to describe the orbit.

If  dynamics are ergodic and lie on a compact space, it can be proved that the entropy of the system  is almost everywhere equal to the orbit complexity. In other words, if such a system  has positive entropy, then for a.e. $x$ the algorithmic information that is necessary to
describe $n$ steps of the orbit of $x$ increases linearly with $n$ and the
proportionality factor is the entropy of the system.
This also implies that if a system has an invariant measure $\mu$, its 
entropy is equal to the mean value of the orbit complexity with respect to $\mu$ \cite{WB}.

%In the following we will consider  dynamical systems $(X,T)$. $X$ is a separable metric space  and $T$ is a function $:X\rightarrow X$. 
It is known  that if a system has
positive entropy the typical speed of separation
of nearby starting orbits is exponential (roughly speaking $\Delta x(t)\simeq\Delta x(0)\ 2^{\lambda t}$).
The speed of separation of exponentially divergent orbits is characterized by the number $\lambda $.
 Lyapunov exponents and Brin-Katok \cite{brinkatok}
local entropy are real valued indicators of the speed
of exponential separation of the orbits (the number $\lambda$ on the exponent of the above formula). The work of Ruelle-Pesin  \cite{pesin}    and
Brin-Katok shows (under some assumption
on the dynamical system) that  their indicators are almost everywhere equal to the
entropy of the system.
In other words: in an ergodic, compact dynamical system the indicator of (exponential) initial data sensitivity is
a.e. equal to the entropy which is a.e. equal to the orbit complexity.
Thus there is a relation between initial data sensitivity and orbit complexity.

This implies that in the compact case orbit complexity and instability are both faces of the same coin. In the case of compact dynamical systems  this motivates  
the general Ford's claim \cite{ford} that orbit complexity was a synonym of chaos.
We remark that if the space is not compact there are examples (\cite{WB}, \cite{io2}) of systems
with no sensitivity to initial conditions  and hight orbit complexity.

 In \cite{io2} another  
definition of orbit complexity is given by the use of computable structures, the definition is given by a different approach involving some
constructivity concept.
This definition is equivalent to Brudno's one if the dynamics lie on a compact space. Replacing compactness with constructivity
allows a more powerful investigation of the relation between orbit complexity
and chaos in the non compact case. 
%In particular, it is shown that
%assuming that the system under consideration is constructive (Definition \ref{definition11} ) we can recover
%a relation between entropy
%and orbit complexity even in the non compact case.

In many  examples of dynamical systems however the entropy could be  0, the speed of
separation of nearby starting trajectories could be less than
exponential and the increasing of the information contained in $n$ step of the orbit could be less than linear. This is the case of the so called {\em Weakly Chaotic Dynamics}. %autoorganizzazione , orlo caos,debolmente caotici.

%first we remark that far to be impossible some different behaviour are%
%discovered to be important in the applications...

%%%%%
The study of weakly chaotic dynamics was discovered to be
important for application purposes.
There are connections with many physical and economic phenomena: self organized criticality, the so called chaos
threshold, the sporadic dynamics, the anomalous diffusion processes  and many others   (see for example \cite{bak},\cite{cont},\cite{tsallis},\cite{Gaspard},\cite{meson}).
In these examples of weakly chaotic dynamics the traditional indicators
of chaos (K.S. entropy, Lyapunov exponents, Brudno's orbit complexity)
vanishes. These indicators are not able to distinguish between all the various
cases of weakly chaotic dynamics. They see all them as trivial dynamical systems.

Some definitions of generalized entropy have been already proposed in 
literature (see e.g. \cite{takens},\cite{meson},\cite{tsallis} for definitions about dynamical systems). One of the most fruitful was the one given by  Tsallis  that found a great
variety of applications: long range interacting  systems, self gravitating
systems, quantum mechanics, social phenomena and many others (see

\noindent http://tsallis.cat.cbpf.br/biblio.htm for an updated bibliography of related topics).
In the literature relations between Tsallis entropy and initial data sensitivity have been  proved (\cite{tsallis},\cite{jin}),
indicating that the main field of application of Tsallis entropy in
dynamical systems is the case of power law initial data sensitivity (if two points starts at distance $\Delta x(0)$ then $\Delta x(t)\sim\Delta x(0) t^{\alpha}$).

%pezzi\#

%pensandoci bene even the Lyapunov exponent...while 
%the K.S. entropy suppposes a linear increasing of the information

%what if the behaviour is not exponential ...

The aim of this work is to give a very general definition of indicators of
initial data sensitivity and orbit complexity and to prove a relation
between them that generalizes the above stated ``orbit complexity=sensitivity \ to initial conditions''.
This will have some interesting corollary as a consequence (Theorem \ref{th41}) and finally we will see some applications of our formulas  to some example of weakly chaotic dynamics.  

 The indicators we will define will have values in 
a totally ordered set ${\cal R   }$ which is constructed in section 3
by the use of the non standard analysis.
In section 2 we add an elementary introduction to non standard analysis, so 
that the paper is self contained.
 The ordered space ${\cal R   }$ will   contain (in some way) a representative of all  the asymptotic  behaviors for  $n\rightarrow \infty$  of  sequences of reals $a_n:{\bf N}\rightarrow {\bf R}$  (for example the various orders of infinity and infinitesimal will correspond to some element of ${\cal R}$). In this way we include in one definition all the possible asymptotic behaviors of initial data sensitivity, orbit complexity or entropy, and so on. 
 
We clarify this by an example. Let us consider the classical definition
of topological entropy for a compact dynamical system $(X,T)$.
We recall and comment the original definition.
If $x,y\in X$ let us say that $x,y$ are $(n,\epsilon)$ separated if  $d(T^k(x),T^k(y))>\epsilon$
for some $k\in\{0,...,n-1\}$. If $d(T^k(x),T^k(y))\leq \epsilon$ for each $k\in\{0,...,n-1\}$ then $x,y$ are said to be $(n,\epsilon)$ near. That is: two 
points are separated if they give rise to substantially different orbits.
A set $E\subset X$ is called $(n,\epsilon)$ separated if $\forall x,y\in E,x\neq y$ then $x,y$ 
 are $(n,\epsilon) $ separated.
 Let us consider
$$ s(n,\epsilon)=max\{card( E):E\subset X \ is \ (n,\epsilon)-separated\}$$

\noindent The number $ s(n,\epsilon)$ measures the number of substantially
different $n$ steps orbits that appear in $(X,T)$.

A chaotic D.S. will have more and more possible different orbits as $n$ increases.
The idea of the definition is that in the more chaotic D.S. (where  the entropy
will be higher)   the cardinality $ s(n,\epsilon)$ increases more quickly as $n$ increases.

The remaining part of the definition is a way to output a real number from this idea.
The number will be a measure of the speed of  exponential increasing of the cardinality
$ s(n,\epsilon)$ as $n$ increases.
We define
$$h(T,\epsilon)=\stackunder{n\rightarrow \infty}{limsup} \frac{\log s(n,\epsilon)}{n}$$

\noindent and then the topological entropy of $T$ is defined as

$$h_{top}(T)=\stackunder{\epsilon \rightarrow 0} {lim}h_K(T,\epsilon).$$

\noindent The logarithm in the definition is taken because $ s(n,\epsilon)$ is expected to increase exponentially $( s(n,\epsilon)\simeq 2^{h(T,\epsilon)n})$.
%In some sense this definition is made ''ad hoc'' for the exponential...
This is a very important case; the case of strongly chaotic dynamic.
However, it is worth remarking that there are examples of chaotic topological dynamical systems with $0$ entropy (see \cite{smital} for example).

When $ s(n,\epsilon)$ increases less than exponentially  the value
of $h(T,\epsilon)$ will be zero.  
If we want to state a definition of generalized topological entropy
that is sensitive to  the cases when $ s(n,\epsilon)$ increases as a power law we could
define (``a la Tsallis'' to some extent)  

$$h_q(T,\epsilon)=\stackunder{n\rightarrow \infty}{limsup} \frac{S_q( s(n,\epsilon))}{n}$$

where $S_q:R\rightarrow R$ is $$S_q(W)=\frac {W^{\frac 1{1-q}}-1}{1-q}.$$

In this   definition the parameter $q$ plays a role similar to the Hausdorff
dimension. Each dynamical system will have a special value of $q$  such that $ S_q( s(n,\epsilon))$ increases linearly and allows a nontrivial value of $h_q$.
The special $q$ of each D.S. will be an indicator of the type of chaotic behavior of the system under consideration.This definition  will classify  the various cases of power law increasing of $ s(n,\epsilon) $  and the exponential one (when $q=1$) 
but it makes no difference between  exponential and stretched exponential\footnote{ $ s(n,\epsilon)\simeq 2^{cn^{\alpha}}.$} (both $q=1$) and between constant and logarithm ($q=0$). 
 
We shall see later with the help of  examples that the asymptotic behavior
of the measures under study exhibit a large variety of different cases. 
So we need the more possible general definition. This will be done with a definition with values in ${\cal R}$.
Moreover in this general setting we can  prove
 very general theorems that will have nontrivial
meaning in all this variety of examples.
We want to remark that the use of the language of the non standard analysis comes about naturally
when we want to consider the  asymptotic behavior  of a sequence of reals.
On the other hand, Benci \cite{bencinsa} gives an elementary but rigorous approach to the nonstandard analysis which we summarize in a  page in section 2.3. This approach is very simple and does not require
any deep tool of logic. So we think that our definitions can be easily understood by readers with no experience in this subject.

In section 3 we define the notions of generalized initial data sensitivity at a point $x$ and  the notion of generalized complexity of the orbit of $x$.
For the definition of constructivity and orbit complexity we introduce the notion of computable structure on a metric space (Section 2.2).
Then, under the assumption that the system is constructive we prove a relation  (Theorem \ref{quellogrosso}) between sensitivity and orbit complexity that is a quantitative 
and rigorous version of the following statement: the asymptotic
behavior of the quantity of information that is necessary to reconstruct
the orbit of $x$ depends on the initial data sensitivity at $x$ and on the
'complexity' of the point $x$ in $X$.
We remark that theorem \ref{quellogrosso} is general, it holds even in infinite
dimensional spaces, provided that the maps we consider are constructive.

This relation will have as a consequence that the points with upper initial 
data sensitivity  ($R(x)$ in Section 3.1) strictly greater than the exponential one are a $0$ dimensional set (Theorem \ref{th41}).
 
We remark that our relations are pointwise. We make no use of invariant measures  that in many cases of weakly chaotic dynamics  are trivial or very
complicated (with multifractal support).
Instead sometime we make use of the natural measure that can be defined
on the metric space $X$: the Hausdorff measure.

In section 4 we apply the the main results of section 3 to give a rigorous
estimation of the orbit complexity in an important example
of weakly chaotic dynamic: the Manneville maps.
 
\section{AIC, computable structures and non standard analysis}

%\centerline{ Ho imparato che tutto il mondo vuole vivere sulla cima della montagna,}
%\centerline{ senza sapere che la vera felicita sta nel modo di scalare la scarpata.}

\subsection{Algorithmic Information Theory}

%We summarize the concepts  leading to the definition of Kolmogorov
%complexity and to the statement of the conjecture (see [2] or [7] for an
%introduction to Algorithmic Information Theory).
%Let $A$ be an algorithm (Partial Recursive Function) that transforms finite
%binary strings to finite binary strings.
%The Kolmogorov complexity $K_A(x)$ of a string $x$ relative to $A$ is
%the length of the shortest string $p$ such that $A(p)=x$, if there are
%no strings $p$ such that $A(p)=x$ then  $K_A(x)=\infty$.
%The string $p$ can be imagined as a program given to a Turing machine
%representing $A$, and the value $A(p)$ can be imagined as the  output
%of the computation.
%If $A$ is an algorithm transforming pairs of binary strings to binary
%strings, the conditional complexity $K_A(x|y)$  of $x$ given $y$ is
%the length of the shortest string $p$ such that $A(p,y)=x$, or if no
%such string exists, it is $\infty$.
%We will consider $K_A(x|n)$, where $n$ is a string representing the binary expan
%sion of     the length of the string
%$x$. In other words $K_A(x|n)$ is the complexity of $x$ given its length.    
%
%{Roughly speaking,a universal algorithm is an algorithm that
%can simulate any other algorithm if an appropriate input is given. For a precise definition see any book of recursion.} {Because if $U$ and $U'$ are un
%iversal then $K_U(s)\leq K_{U'}(s)+c$ where $c$ is a constant depending only on
%$U$ and $U'$, this tells that when using the universal algorithm $U$ the complexity of $s$ with respect to $U$ depends only on $s$ up to a fixed constant.}.

In this section we give an introduction to algorithmic information theory.
The introduction will be informal, to help the reader that is not 
familiar with recursion theory to understand the paper.
A more detailed exposition of algorithmic information theory can be found in \cite{Zv} or \cite{Ch}. 

Let us consider the set  $\Sigma=\{0,1\}^*$ of finite (possibly empty)  binary strings. If $s$ is a string we define $|s|$ as the length of $s$.

Let us consider a Turing machine (a computer) $C$: by writing $C(p)=s$ we mean that $C$ starting with input $p$ (the program) stops with output $s$ ($C$ defines a  partial recursive function  $C:\Sigma \rightarrow \Sigma$). If the input gives a never ending computation the output (the value of the recursive function) is not defined.
If $C:\Sigma \rightarrow \Sigma$ is recursive  and its value is defined
for all the input strings in $\Sigma$ (the computation stops for each input)
then we say that $C$ is a {\em total} recursive function from $\Sigma$ to $\Sigma$. The algorithmic information content of a string will be the length of the shortest
program that outputs the string.

\begin{definition}
The Kolmogorov complexity or algorithmic information content of a string 
$s$ given $C$ is the length of the smallest
program $p$ giving $s$ as the output:

$$
{ K}_C(s)=\stackunder{C(p)=s}{\min }|p| 
$$

\noindent  if $s$ is not a possible output for
the computer $C$ then ${K}_C(s)=\infty $ .
\end{definition}
For example by this definition we see that the algorithmic information content (A.I.C.)
of a $2n$ bits   long periodic string 
$$s=''1010101010101010101010...''$$  is small because 
the string is output of a shortest program: 

\centerline{\em repeat $n$ times (write (``10''))}

\noindent the AIC of the string $s$ is then less or equal than $ log(n)+Constant$ where $log (n)$ bits are sufficient  to  code ``$n$''  and the constant represents the length of the code for the computer $C$ representing the instructions ``repeat...''. As it is intuitive the information content of a periodic string is very poor.
On the other hand each $n$ bits long string $$s'=''1010110101010010110...''$$ is output of the trivial program 
$${write(''1010110101010010110...'')}$$ 
this is of length $n+constant$ this implies that the AIC of each string
is (modulo a constant which depends on the chosen computer $C$) less or equal than its length.

Until this point the algorithmic information content of a string depends on the choice of the computer $C$. We will see that there is a class of computers that 
allows an ``almost'' universal definition of algorithmic information content of a string: if we consider computers from this class the A.I.C. of a string
will be defined independently of the computer up to a constant. In order
to define such a class of universal computers 
 we give some notations that are necessary to work with strings:
there is a correspondence  $c:\Sigma \rightarrow {\bf N}$ from the set $\Sigma$  and the set \-${\bf N}$ of natural numbers 

$$ \emptyset \rightarrow 0, \ \ \  0 \rightarrow 1 , \ \ \  1
\rightarrow 2, \ \ \  00 \rightarrow
3, \ \ \  01\rightarrow 4 , \ \ \ 10\rightarrow 5 , ...$$

 This correspondence  allows us to interpret  natural numbers as  strings
and vice-versa when it is needed. 
 We remark  that
$|s|\leq log(c(s))+1$. \footnote{ In this paper all the logarithms are in base two.}

If  $s$ is a string with $|s|=n$ we denote by  $s\hat{\ \ } $ the string $s_0s_0s_1s_1...s_{n-1}s_{n-1}01.$ If $a=a_1...a_n$ and $b=b_1... b_m$ are strings
then $ab$ is defined as the string $a_1...a_nb_1... b_m$. If $a$ and $b$ are strings
then $a\hat{\ \ } b$ is an encoding of the couple $(a,b)$. There is
an algorithm that getting the string $a\hat{\ \ } b$ is able to recover both the strings $a$ and $b$. 
An universal Turing machine intuitively is a machine that can emulate 
any other  Turing machine if an appropriate input is given. 
We recall that there is a recursive enumeration $A_1,A_2...$
of all the Turing machines.    
\begin{definition}
A  Turing machine ${\cal U}$ is said to be universal if for all $m\in {\bf N}$ and $p\in \Sigma$ then  ${\cal U}(c^{-1}(m)\hat{\ \  } p)=A_m(p)$.

\end{definition}

%vecchia def.
%\begin{definition}
%A  Turing machine ${\cal U}$ is said to be universal if for each
% function ${\cal F}$ there exist a prefix $\mu $ such that
%for each program $p$: ${\cal U}(\mu p) ={\cal F}(p)$.
%\end{definition}
In the last definition the machine ${\cal U}$ is universal because ${\cal U}$ is able to emulate each other machine
$A_m$ when in its input we specify the number $m$ identificating $A_m$ and
the program to be runned by $A_m$.
It can be proved  that an universal
Turing machine 
exists.

\begin{definition}
A  Turing machine $ F$ is said to be {\em asymptotically optimal }if for each Turing machine $C$ and each binary string $s$ we have $K_{F}(s)
\leq {K}_C(s)+c $ where
the constant $c$ depends on $C$ and not on $s$.
\end{definition}

The following proposition can be proved from the definitions

\begin{proposition}
If ${\cal U}$ is an universal Turing machine
then ${\cal U}$ is asimptotically optimal.
\end{proposition}
 
 This tells us that choosing an universal Turing machine the
complexity of a string is defined independently of the given 
Turing machine  up to a constant.  
 For the remaining part of the
paper we will suppose that  an universal Turing machine ${\cal U}(p)$ is chosen once forever.

\subsection{Computable  Structures, Constructivity }

A computable structure on a separable metric space $(X,d)$ is a class of dense immersions ($I:\Sigma\rightarrow X$)  of the
space of finite strings $\Sigma$ in the metric space. The immersions are such that the distance
$d$ restricted to the points that are images of strings ($x=I(s):x\in X,s\in\Sigma$) is a ``computable'' function.
Many concrete metric spaces used in analysis or in geometry have a
natural choice of a computable structure.
The use of computable structures allows to consider algorithms acting
over metric spaces and to define constructive functions between metric spaces,
that is, functions such that we can work with by using a finite amount 
of information. In the following we often will assume that the
dynamical systems under our consideration are constructive. All the dynamical system that we can construct explicitely 
are construcive.
From the philosophical point of wiew we  think  that the assumption of constructivity is not unnatural because
 even if the  maps coming from  physical reality were not constructive, the models used  to describe such a reality should be constructive (to allow calculations).
On the other hand, to add constructivity  allows to prove stronger
theorems, avoiding pathologies coming from random maps.

An interpretation function is a way to interpret  a string as a point of the metric space.

\begin{definition}
An interpretation function on  $(X,d)$ is a function $I:\Sigma \rightarrow
X$ such that $I(\Sigma)$ is dense in $X$.
\end{definition}

 A point $x\in X$ is said to be {\em ideal} if it is the image of some string $x=I(s), s\in \Sigma$.  
 An interpretation is said to be computable  if the distance between ideal points is computable with arbitrary precision:

\begin{definition}
A computable  interpretation function on  $(X,d)$ is a function $I:\Sigma \rightarrow
X$ such that $I(\Sigma)$ is dense in $X$ and  there exists a total recursive
function $D:\Sigma \times \Sigma \times
 {\bf N} \rightarrow {\bf Q}$ such that $\forall s_1,s_2 \in \Sigma ,n\in
{\bf N}$:

$$|d( I(s_1),I(s_2))-D(s_1,s_2,n)|\leq \frac{1}{2^n}.$$

\end{definition}

Two interpretations are said to be equivalent if the distance from an ideal
point from the first and a point from the second is computable up to
arbitrary precision.
For example, the finite binary  strings $s\in \Sigma $ can be interpreted as rational numbers by interpreting the string as the binary expansion of a number. Another interpretation can be given by
interpreting a string as an encoding of a couple of integers whose ratio
gives the rational number. If the last encoding is recursive, the two interpretation are equivalent.

\begin{definition}\label{def7}     
Let $I_1$ and $I_2$ be two computable interpretations in $(X,d)$;
we say that $I_1$ and $I_2$ are equivalent if there exists a total recursive function
$D^*:\Sigma \times \Sigma \times {\bf N} \rightarrow {\bf Q}$,
such that $\forall s_1,s_2 \in \Sigma ,n{\in {\bf N}}$:

$$|d( I_1(s_1),I_2(s_2))-D^*(s_1,s_2,n)|\leq \frac{1}{2^n}.$$

\end{definition}

\begin{proposition}
The relation  defined by definition \ref{def7} is an equivalence relation.
\end{proposition}

\noindent For the proof of this proposition see \cite{io2}.

\begin{definition}
A computable structure ${\cal I}$ on $X$ is an equivalence class  of computable  interpretations  in $X$.
\end{definition}

%$I(s)=s(1)...s([n/2]).s([n/2]+1)...s(n).$
For example if $X={\bf R}$ we can consider the interpretation $I:\Sigma \rightarrow {\bf R}$ defined in the following way: if $s=s_1...s_n\in \Sigma$ then \begin{equation}\label{I}I(s)=\sum_{1\leq i \leq n}s_i2^{[n/2]-i}.\end{equation} This is an interpretation of a string
as a binary expansion of a number.
  $I$ is a computable interpretation, the computable structure on ${\bf R}$ containing $I$ will be called {\em standard} computable structure.
If $r=r_1r_2...$ is an infinite string such that $lim \frac{K_C(r_1...r_n)}{n}=1$\footnote{such a string exist, see for example \cite{io1} theorem 13.} then the interpretation $I_r$ defined as $I_r(s)=I(s)+\sum r(i)2^{-i}$ is computable but not equivalent to $I$. $I$ and $I_r$ belongs to different computable structures.
 
In a similar way it is easy to construct computable structures in ${\bf R}^n$ or in separable function spaces codifying a dense subset (for example the set of step functions) with finite strings. 
We remark as a property of the computable structures that if $B_r(I(s))$ is 
an open ball with center in an ideal  point $I(s)$ and rational radius $r$ and $I(t)$ is another point   then there is an algorithm that 
verifies if  $I(t)\in B_r(I(s)) $. If  $I(t)\in B_r(I(s)) $ then the algorithm
outputs ``yes'', if $I(t)\notin B_r(I(s)) $ the algorithm  outputs ``no'' or 
does not stop. The algorithm calculates $ D(s,t,n)$ for each $n$ until it
finds that $D(s,t,n)+2^{-n}< r$ or  $D(s,t,n)-2^{-n}> r$, in the first case
it outputs ``yes'' and in the second it outputs ``no'', if $d(I(s),I(t))\neq r$ the algorithm will stop and output an answer.

We give a definition of {\em morphism} of metric spaces with computable 
structures, a morphism is heuristically  a computable
function between  computable metric spaces.

\begin{definition}
If $(X,d,{\cal I})$ and $(Y,d',{\cal J}) $ are spaces with  computable
structures; a function $\Psi :X\rightarrow Y$ is said to be  a morphism of
computable structures if $\Psi $ is uniformly continuous and for each pair
$I\in {\cal I},J\in {\cal J}$ 
there exists a total recursive function
$D^*:\Sigma \times \Sigma \times {\bf N} \rightarrow {\bf Q}$,
such that $\forall s_1,s_2 \in \Sigma ,n{\in {\bf N}}$:

$$|d'(\Psi( I(s_1)),J(s_2))-D^*(s_1,s_2,n)|\leq \frac{1}{2^n}.$$
%$\Psi  (I(*)):\Sigma \rightarrow Y$
%is an interpretation function    equivalent to $J$.
\end{definition}

\noindent We remark that $\Psi $ is not required to have dense image
and then $\Psi(I(*))$ is not necessarily an interpretation function 
equivalent to $J$.

\begin{remark}\label{remark10} 
As an example of the  properties of the morphisms, we remark that if a map $\Psi:X\rightarrow Y $  is a morphism  then given a point $x\in I(\Sigma) \subset X$ it is possible to find by an algorithm a point $y\in J(\Sigma) \subset Y$ as near as we want to $\Psi (x) $. 
\end{remark}
The procedure is simple: if $x=I(s)$ and we want to find a point $y=J(z_0)$ such that $d'(\Psi( I(s)),y)\leq 2^{-m} $ then we calculate $D^*(s,z,m+2)$ for each $z\in \Sigma$ until we find $z_0$ such that $D^*(s,z_0,m+2)<2^{-m-1}$. Clearly $y=J(z_0)$ is such that $d'(\Psi(x),y)\leq 2^{-m} $. The existence of such a $z_0$ is assured by the density of $J$ in $Y$.
In particular the identity is a morphism. Remark \ref{remark10}  applied to the identity  will be used in the proof of lemma \ref{lemma12}.

%metti traslazioni irrazionali

A constructive map is a morphism for which the continuity relation between $\epsilon$ and $\delta$  is given by a
 recursive function. The following is in some sense a generalization of the definition of Grzegorczyk, Lacombe (see e.g. \cite{purel}) of constructive function.

\begin{definition}\label{definition11}
A function $\Psi:(X,d,{\cal I})\rightarrow (Y,d',{\cal J}) $ between  spaces with computable structure
$ (X,{\cal I}) ,(Y,{\cal J})$ is said to be constructive
if  $\Psi$ is a morphism between the  computable structures and it is 
effectively uniformly continuous, i.e. there is a total recursive function
$f:{\bf N}\rightarrow {\bf N} $ such that for all $x,y \in X$  
$d(x,y)<2^{-f(n)}$ implies $d'(\Psi (x),\Psi (y))<2^{-n}$.
\end{definition}

If $f(n)$ is recursive and satisfies the hypothesis that, $d(x,y)<2^{-f(n)}$ implies $d(T (x),T (y))<2^{-n}$ then  $max (f(n),n)$ is recursive and still satisfies the hypothesis. Then we can suppose that if $f(n)$ is a function of effective continuity
then $f(n)\geq n$.

If a map between spaces with a computable structure is constructive then there
is an algorithm to follow the orbit  each ideal point $x=I(s_0)$.

\begin{lemma}\label{lemma12}
If $T: (X,{\cal I})\rightarrow (X,{\cal I})$  is constructive, $ I\in{\cal I}$ then  there is an algorithm (a total recursive function)
$A:\Sigma \times {\bf N}\times {\bf N}\rightarrow \Sigma$ such that $\forall k,m\in {\bf N},s_0\in \Sigma$ $d(T^k(I(s_0)),I(A(s_0,k,m)))< 2^{-m} $. 
\end{lemma}

Proof. 
Since $T$ is  effectively uniformly continuous
we define the function $g_k(m)$ inductively as 
$g_1(m)=f(m )+1$, $g_i(m)=f(g_{i-1}(m)+1)$ where $f$ is the function of effective uniform
continuity of $T$ (definition \ref{definition11}).
 If  $d(x,y)<2^{-g_k(m)} $ then $d(T^i(y),T^i(x))<2^{-m}$ for $i\in \{1,...,k\}$. 
Let us choose $I \in {\cal I}$. 
We recall that the assumption that $T$ is a morphism implies that
there is a recursive function $D^* (s_1,s_2,n)$ such that 

$$|D^* (s_1,s_2,n)-d(I(s_1),T(I(s_2)))| <2^{-n}.$$

\noindent Let us suppose that $x=I(s_0)$.
Now let us describe the algorithm $A$:
 using the function  $D^*$ and the function $f$,  $A$ calculates $g_k(m)$ and finds a string $s_1$ such that 
$d(I(s_1),T(I(s_0))) <2^{-g_k(m)}$ as described in remark \ref{remark10}. This is the first step of the
algorithm. Now $d(T(I(s_1)),T^2(x))\leq 2^{-(g_{k-1}+1)}$.
We can use $D^*$ to find a string $s_2$ such that $d(I(s_2),T(I(s_1)))<2^{-(g_{k-1}+1)}$. By this $d(I(s_2),T^2(x)) \leq 2^{-g_{k-1}}$. 
 This implies that $d(T(I(s_2)),T^3(x))\leq 2^{-(g_{k-2} +1)}$, then 
we find $s_3$ such that $d(I(s_3),T(s_2))\leq 2^{-(g_{k-2} +1)}$ and so on for
$k$ steps. At the end we find a string $s_k$ such that $d(I(s_k),T^k(x)) \leq 2^{-m}$. $\Box $

\subsection{Non standard analysis}

We define the extended real line ${\bf R}^*$ to be an ordered field satisfying suitable axioms. The existence of such a field is proved in \cite{bencinsa}.
 ${\bf R}^*$ will contain  the standard real numbers and other elements representing 
  {\em infinite} and the {\em infinitesimal} numbers.

We call Hyperreal Line a field satisfying the following axioms:
\begin{axiom}  The set of the hyperreal numbers ${\bf R}^{*}$ is an ordered field which
contains ${\bf R}$ as a subfield.
\end{axiom}

\begin{axiom}\label{ax15} There is a surjective ring homomorphism

\begin{center}
$J:{\bf R}^{\bf N}\rightarrow {\bf R}^{*}$
\end{center}associating to each real sequence an hyperreal number.

\end{axiom} 

Intuitively the homomorphism $J$ associates to a sequence of reals its asymptotic behavior. For example if $\stackunder{n\rightarrow\infty}{lim} a_n=0$ then $J(a_n)$ will be an {\em infinitesimal} number. Moreover if $a=(a_i)$ and $b=(b_i)$ are two
sequences we would like that if $a_i\geq b_i$ for all $i \in {\bf N}$ then $J(a)\geq J(b)$.
For this reason  $J$ is required to satisfy the following monotonicity property
 
\begin{axiom}\label{ax16}
 If there exists $k\in {\bf N}^{+}$ such that 

\begin{center}
$\forall n\in {\bf N}^{+},\phi _{kn}\geq a$
\end{center}

\noindent with $a\in {\bf R}$, then 

\begin{center}
$J(\phi )\geq a$ . 
\end{center}

\end{axiom}
 From axiom \ref{ax16} it follows for example that
 if $a_i=2^{-i}$ and $b_i=\frac 1i$ then $0\leq J(a_i)\leq J(b_i)$.
Another consequence of axiom  \ref{ax16} is that  for each real sequence $x_i$ $liminf (x_i)\leq J(x_i) \leq limsup(x_i)$.

A field satisfying our axioms exists. As the reader could imagine, it can be constructed from  the set of real sequences modulo a suitable equivalence relation (\cite{bencinsa}).
%Now we introduce a special number $\alpha _0$ if $i:N\rightarrow R$ is the
%natural imbedding $i(n)=n$ then $\alpha _0=J(i)$
%
%Dire che ogni numero si vede comeuna funzione di $\alpha _0$
%

\subsubsection{Extension of functions, infinite and infinitesimal numbers.}

%Given any function $\phi :N\rightarrow R$ we extend it to a function $\phi
%:N\cup (\alpha _0)\rightarrow R^{*}$
%
%\begin{center}
%$\phi (\alpha _0)=J(\phi ).$
%\end{center}
%
%natural preextension of $\phi .$

Given any function  $f:{\bf R}\rightarrow {\bf R}$  we extend it to a function  $%
f^{*}:{\bf R}^{*}\rightarrow {\bf R}^{*}$ as follows: if $a_i$ is a real sequence
and $x=J(a_i)\in {\bf R}^*$ we define 
\begin{center}
%$f^{*}(\phi (\alpha _0))=J(\phi (\alpha _0)).$
$f^*(x)=f^{*}(J(a_i))=J(f (a_i)).$
\end{center}

\begin{proposition}\label{hk}
 the definition is well posed i.e. % $\phi (\alpha _0)=\psi (\alpha
%_0)$ implies that 
$J(a_i)=J(b_i)$ implies that $J(f(a_i))=J(f(b_i)).$
\end{proposition}
The proof of proposition \ref{hk} can be found in \cite{bencinsa}.
%\begin{center}
%$f^{*}(\phi (a_0))=$ $f^{*}(\psi (a_0))$
%\end{center}

As we stated before, in ${\bf R}^*$ there are some elements representing 
 the infinite and infinitesimal numbers:

\begin{definition}
An hyperreal number $\xi $ is called infinite if $\forall k\in {\bf N}$ we have $|\xi |>k$. A number $\xi $ is called infinitesimal if $\forall k\in {\bf N},|\xi |<\frac 1k$.
A number $\xi $ is called bounded if $\exists k\in {\bf N},|\xi |<k$.

%$\alpha _0$ is infinite, $\frac 1{\alpha _0}$ is infinitesimal.
\end{definition}
For example the reader could verify directly from the axioms that if $a_i:{\bf N}\rightarrow {\bf R}$ is the identity: $a_i=i$ then $J(a_i)$ is infinite and $J(\frac 1 {a_i}) $ is infinitesimal.

\section{sensitivity and orbit complexity}

Now we construct the space $\cal R$ in which our indicators of orbit complexity
and initial data sensitivity will have value.

\begin{definition} 
If $a,b\in {\bf R}^*$ we say that $a$ and $b$ have the same order 
and write $a\simeq b$ if  and only if both $\frac ab$ and $\frac ba$ are bounded.
\end{definition}
 $\simeq $ it is clearly an equivalence relation. In the following  by $[\alpha ]$  we will indicate the equivalence class of $\alpha$. 

We now define an ordering relation on the quotient space  $\frac {{\bf R}^*}{\simeq}$.
We say that $[a]\leq [b]$\footnote{By an abuse of notations we use the symbol $\leq$ for this ordering relation, this will cause no ambiguity with the ordering relation defined on ${\bf R}^*$. } if $\forall x\in [a] ,\forall y\in [b]$ then or $x\simeq y$ or $x< y$. We remark that the order relation on ${\bf R}^*$ is compatible with
the equivalence relation $\simeq$, if $a\leq b $ in ${\bf R}^*$ then $[a]\leq [b]$ in $\frac {{\bf R}^*}{\simeq}$ .  And thus  $\frac {{\bf R}^*}{\simeq}$
is totally ordered. The relation $<$ is then  defined in the obviuos way as $[a]<[b]\iff [a]\leq[b],[a]\neq [b]$.

$\frac {{\bf R}^*}{\simeq}$ contains a representative of all the infinite (infinitesimal) asymptotic behaviors of real sequences. $\frac {{\bf R}^*}{\simeq}$ is sometime called the group of orders.
The natural projection from ${\bf R}^*$ to $\frac {{\bf R}^*}{\simeq}$ $:a\rightarrow [a]$ allows to forget
all the lower order terms in the  hyperreal number $a$: for example if $a_i=i$ as above $J({a_i}^2+a_i)$ and $J({a_i}^2+\sqrt{a_i}) $ belongs to the same class as $J({a_i}^2)$, in other words $[J({a_i}^2+a_i)]=[J({a_i}^2+\sqrt{a_i})]=[J({a_i}^2)]$.

Unfortunately $\frac {{\bf R}^*}{\simeq}$ is not complete (as ${\bf R}^*$ is not complete). The supremum or the infimum of a sequence in $\frac {{\bf R}^*}{\simeq}$ may not exist in $\frac {{\bf R}^*}{\simeq}$.
For this reason we will consider   a space $\cal R$ which is a completion of 
$\frac {{\bf R}^*}{\simeq}$,
the {\em sup} and {\em inf} of each sequence in $\frac {{\bf R}^*}{\simeq}$
is in $\cal R$.
We now outline a possible construction of a completition of $\frac {{\bf R}^*}{\simeq}$, there are other possible costructions. Another possible completition of $\frac {{\bf R}^*}{\simeq}$ can be constructed for example by Dedekind sections.
We construct $\cal R$ by quotienting the set of monotone sequences in $\frac {{\bf R}^*}{\simeq}$ by a suitable equivalence relation. This will add the supremum to each countable sequence.
\begin{proposition}
There is an ordered space ${\cal R}$ such that $\frac{{\bf R}^*}\simeq\subset{\cal R}$ in a natural and order preserving  way and if $a_i\in \frac{{\bf R}^*}\simeq $ is a monotone sequence then $inf ( a_i)$ and $sup (a_i)$ are in ${\cal R}$. 
\end{proposition}
{\em Proof.}
Let us consider the set of monotone sequences in $\frac{({\bf R}^*)}\simeq$: $${\cal A}=\{(a_i):{\bf N }\rightarrow {\frac{{\bf R}^*}\simeq} \ s.t. \ a_i \ is \  monotone\} $$ in $\cal A$ the set  $\frac{{\bf R}^*}\simeq$ will be identified with the subset 
of constant sequences.

We define the ordering relation on  ${\cal A}$ in the following way
$(a_i)< (b_j)$ $\iff$ $\exists M \ $ such that $ \ \forall n,m  $ with $n>M, m>M $ then $a_n< b_m$ .

We define the relation $\approx$ in the following way: $(a_i) \approx  (b_j)$
if  neither $(a_i)< (b_j)$ nor $(b_j)< (a_i)$. $\approx$ is an equivalence relation:  
$\approx $ is trivially symmetric and reflexive. The transitivity follows
from the remark that if $a_i$ is not $<b_j $ then $\forall M\ \exists n,m  \ s.t. \ n>M, m>M , a_n \geq b_m $.
As it is easy to verify this is a transitive relation (because the sequences in $\cal A$ are monotone).

The set of equivalence classes ${\cal R}=\frac{\cal A}{\approx}$ is then totally
ordered in the same way as before and contains ${\frac{{\bf R}^*}\simeq}$ as a subset. Moreover 
each monotone sequence in ${\frac{{\bf R}^*}\simeq}$ has its $sup$ and $inf$
in ${\cal R}$: let us indicate with $p$
the {\em natural projection} associating to each element in ${\cal A}$ its equivalence class.
 If for example  $a=(a_i)$ is a nondecreasing sequence then $inf(a) $ is the equivalence class of the constant sequence $b=(b_j)$ such that $\forall j\ b_j=a_0$ and $sup(a)=p(a)$  ($p$ is the natural projection map as defined above). It is easy to verify that
$p(a)>a_j$ for each $ j\in {\bf N}$  and for each $ x \in {\cal R}\ s.t. $ $\forall j\in {\bf N}, x> a_j \ $ we have $x> p(a)$ or $x\approx p(a)$ .$\Box$

The set $\cal R$ may look in some way mysterious and not practical to be used.
The reader will see in the examples that the elements of $\cal R$ we will have as value of our invariants will be classes that can be expressed in an explicit
way. For example a possible value of $r(x)$ (definition \ref{d23})  could be $[J(n^{-\frac 12})]$ (the class in $\frac{R^*}\simeq$ containing the asymptotic behavior  of the sequence $a_n=n^{-\frac 12}$). Since $\frac{R^*}\simeq$ is immersed in ${\cal R}$ in a natural way we can consider an element of $\frac{R^*}\simeq$ as an element of $\cal R$ without ambiguity.  Another possible value of $r(x)$ could be  $[J(n^{-\frac 13})]$ and it is clear
that $[J(n^{-\frac 12})]\leq[J(n^{-\frac 13})]$ so we can easily compare the values.

The notion of infinite and infinitesimal numbers can be extended to 
the elements of $\cal R$:
if $v\in {\bf R},v>0$ then  $c=[J(v)]$  is the element of ${\cal R}$
corresponding to the class of bounded and not infinitesimal numbers.
Moreover:
an element  $\epsilon \in {\cal R}  $ is said to be infinitesimal if $\epsilon<c$
and
an element  $\gamma \in {\cal R}  $ is said to be infinite if $\epsilon>c$.

Finally we remark that ${\cal R}$ is closed by countable  $inf$ and $sup$
and the projection $p$ (defined in the proof of the above proposition) can be extended to a function from the set of monotone sequences
in $\cal R$ to $\cal R$, associating to a sequence its supremum or infimum
according that the sequence is increasing or decreasing.

\begin{proposition}
If $(a_i)$ is a monotone sequence in $\cal R$, then $$sup (a_i), inf(a_i)\in {\cal R}$$.
\end{proposition}
{\em Proof.} Let us suppose that $(a_i)$ is a non decreasing sequence, the case where $(a_i)$ is non increasing is analogous. We show that $sup (a_i)\in {\cal R}$. If $(a_i)$ is eventually constant the proposition is obvious.
If $(a_i)$ is not eventually constant let us consider a subsequence $a_{i_k}$
that is strictly increasing ($\forall k$ $ a_{i_k}>a_{i_{k-1}}$).
The elements $ a_{i_k},a_{i_{k-1}}\in {\cal R} $ are equivalence classes of monotone sequences with values in $\frac {R^*}{\simeq}$. Let use consider two of this sequences $(\alpha ^{k-1}_n)_{n\in{\bf N}}:{\bf N}\rightarrow \frac {R^*}{\simeq}$ such that $p((\alpha ^{k-1}_n)_{n\in{\bf N}})= (a_{i_{k-1}}) $
and  $(\alpha ^{k}_n)_{n\in{\bf N}}:{\bf N}\rightarrow \frac {R^*}{\simeq}$ such that $p((\alpha ^{k}_n)_{n\in{\bf N}})= ( a_{i_{k-1}}) $. Since $(\alpha ^{k-1}_n) <(\alpha ^{k}_n)$ then $\exists M_k\in {\bf N} \ s.t. \ \forall n \  \alpha ^{k-1}_n <\alpha ^{k}_{M_k} $.  This is true for each choice of the sequences $\alpha ^k,\alpha ^{k-1}$ in the classes $a_k,a_{k-1}$, hence what follows does not depends on the choice of $\alpha ^k$ and $\alpha ^{k-1}$. Let us now consider the sequence $(\alpha ^{k} _{M_k})_{k\in{\bf N}}:{\bf N}\rightarrow \frac {R^*}{\simeq}$ this is an increasing sequence and $\forall j<k , \ a_{i_j}<\alpha ^k _{M_k}$, then $\forall i\in {\bf N}, p(\alpha ^{k} _{M_k})>a_i $. Now it is easy to see that if  $ b $ is such that
$\forall i_k b>a_{i_k}$ then $b\approx p(\alpha ^k _{M_k})$ or $b>p(\alpha ^k _{M_k})$, then $ p(\alpha ^k _{M_k})$ is  $sup (a_i)$. $\Box$  

By the above result we also see that the function $p:{\cal A}\rightarrow {\cal R}$ can be extended to a function $\overline {p}:\{(a_i):{\bf N }\rightarrow {\cal R} \ s.t. \ a_i \ is \  monotone\}\rightarrow {\cal R} $ by associating to each sequence its $sup$ of $inf$ according that the sequence is increasing or decreasing.

\subsection{Initial data sensitivity}

Let $X$ be a separable metric space  and $T$ a function $X\rightarrow X$.
\footnote{We remark that in this and in the following subsection we do not require that $T$ is continuous.}
Let us consider the following set:

$$ B(n,x,\epsilon)=\{ y \in X :d(T^i(y),T^i(x))\leq \epsilon \ \forall i \ s.t. \ \ 0\leq i \leq n\}.$$

$ B(n,x,\epsilon)$ is the set of points ``following'' the orbit of $x$ for $n$
steps at a distance less than $\epsilon$. When the orbits of $(X,T)$ diverges  the set $ B(n,x,\epsilon)$ will be smaller and  smaller as $n$ increases. 
 The speed of decreasing of the size of this set considered
as a function of $n$ will be a measure of the sensitivity of the system
to changes on initial conditions.

Brin and Katok used the set $ B(n,x,\epsilon)$ for their definition 
of local entropy \cite{brinkatok}. In their paper the measure of the size of $ B(n,x,\epsilon)$ 
was the invariant measure of the set.

If we are interested to approximate the orbit of $x$ for $n$ steps we are interested to know how close we must approach the initial condition $x$ 
to ensure that the resulting approximate orbit is close to the orbit of $x$; another possible measure of the size of $ B(n,x,\epsilon)$ is then the radius of
the biggest ball with center $x$ contained in $ B(n,x,\epsilon)$.

$$ r(x,n,\epsilon)=\stackunder {B_r(x)\subset B(n,x,\epsilon)}{ sup\ r}.$$

Or the radius of the smaller ball that contains $ B(n,x,\epsilon)$ 

$$ R(x,n,\epsilon)=\stackunder {B_R(x)\supset B(n,x,\epsilon)}{ inf\ R}.$$

As said before the generalized initial data sensitivity will be a function associating to a point of $X$ a class in ${\cal R}$ indicating how faster orbits coming  from a neighborhood of $x$ will diverge. For this purpose we measure how faster
 $r(x,n,\epsilon )$ decreases as $n$ increases, i.e. we consider the asymptotic behavior of the sequence  $r(x,n,\epsilon )$ as $n$ increases.
First we define

\begin{definition}We define
$r_{\epsilon}:X \rightarrow \frac {R^*}{\simeq} $ as

$$r_\epsilon (x)=[J(r(x,n,\epsilon))]$$

\noindent and

$R_{\epsilon}:X \rightarrow \frac {R^*}{\simeq} $ as

$$R_\epsilon (x)=[J(R(x,n,\epsilon))]$$

\end{definition}

The following lemma implies that $r_\epsilon(x) $ and $R_\epsilon(x) $ are  monotone functions with respect to $\epsilon$.

\begin{lemma}
If $\epsilon>\theta$ then $r_{\epsilon}(x) \geq r_{\theta}(x)$, $R_{\epsilon}(x) \geq R_{\theta}(x)$.
\end{lemma}

Proof. Obvious $\Box$

For the previous lemma we define the indicator of initial data sensitivity
at $x $ by letting $\epsilon $ go to $0$ as the infimum of the $r_{\epsilon}(x) $  for $\epsilon\in {\bf R},\epsilon >0$. 
This infimum will be in $\cal R$.

%Unfortunately we cannot assure
%that in principle the infimum of a set of hyperreal numbers is in ${\bf R}^*$.
 %So we consedered the set ${\cal R} $ that is constructed in a way that
%it contains all the infimums of sequences in $\frac{R^*}{\simeq}$.
%
% the set ${\cal R}/\simeq$ must be imagined as the set of all the {\em different} possible way to go 
%at 0 or $\infty$ by a sequence of reals
%

\begin{definition}\label{d23}

We define the indicator of initial data sensitivity at $x$ as

$ r(x):X\rightarrow {\cal R}$

$$r(x)=\stackunder {\epsilon \in {\bf R}^+}{inf}r_{\epsilon}(x)$$

\noindent In the same way we define

$ R(x):X\rightarrow {\cal R}$
$$R(x)=\stackunder {\epsilon  \in {\bf R}^+ }{inf}R_{\epsilon}(x).$$

\end{definition}

The classical definition of dynamical system sensitive to initial 
conditions is related to our last  definition. To say that a system is sensitive to initial conditions is equivalent to say that there is a $\delta$ such that  $r_{\delta}(x)$ is infinitesimal for all
the $x\in X$:

\begin{definition}
A dynamical system $(X,T)$ is said to have sensitive dependence on initial conditions
if there is a $\delta$ such that for each $x\in X$ and every neighborhood $U$ of $x$ there is $y\in U$ and $k\in {\bf N}$ such that $d(T^k(x),T^k(y))>\delta$. 
\end{definition}

 \begin{proposition}
A system has sensitive dependence on initial conditions if and only if
there is a $\delta$ such that $\forall x \in X$ $r_{\delta}(x)$ is infinitesimal.
\end{proposition}
\noindent The proof follows directly from the definition of $r_{\delta}(x)$.
%By lemma 14 we know that $r(x)=[\stackunder {\epsilon \rightarrow 0}{lim} r_{\epsilon}(x)]$ and  $R(x)=[\stackunder {\epsilon \rightarrow 0}{lim} R_{\epsilon}(x)]$.
%In appendix 2 we outline the construction that is necessary to enlarge $R^*$ 
%and obtain a definition of $r(x) $ and $R(x)$ that is well defined on the
%whole space $X$.

We give some example of different behaviors of $r(x)$ and $R(x)$ in
dynamical system over the interval $[0,1].$
\noindent {\em The identity map} $T(x)=x$. In this map $\forall n$ $B(n,x,\epsilon )=\{y\in
[0,1],s.t.|y-x|<\epsilon \}$ then if we choose for example $x=\frac 12$ we
have $R(\frac 12,n,\epsilon )=r(\frac 12,n,\epsilon )=\epsilon $ i.e. the constant sequence with
value $\epsilon .$ Then $R_\epsilon(\frac 12)=r_\epsilon (\frac 12)=[J(\epsilon )]$ where $J(\epsilon)$ is the number in $R^{*}$ corresponding to the constant sequence with value $\epsilon$  and $R(\frac 12)=r(\frac 12)=[J(\epsilon )]=\{x\in
R^{*}\ s.t. \ x\ is\ bounded\}$ i.e. the class containing the numbers in $R^{*}$
corresponding to the constant sequences$.$ The same arguments can be applied
to the irrational translation on $[0,1]$: $T(x)=x+t \ (mod  \ 1 )$ where $t\notin {\bf Q}$ obtaining the same kind of initial data sensitivity as the identity (in effect both the maps are not sensitive to initial conditions).   

\noindent {\em The one dimensional baker's map} $T:[0,1]\rightarrow [0,1],T(x)=2x\ (mod\ 1).$ If we choose  for example $x=0$ , we
have  $B(n,0,\epsilon )=\{y\in [0,1], 0 \leq y \leq 2^{-n}\epsilon \},$ $R_\epsilon(0)=r_\epsilon (0)=[J(\epsilon
2^{-n})]$ and $R(0)=r(0)=[J(2^{-n})]$ i.e. the class containing all the bounded
multiples of the exponential infinitesimal number.

\noindent {\em The piecewise linear map}  $T:[0,1]\rightarrow [0,1]$ 

$T(x)=\left\{ \begin{array}{cc}
\frac{\xi _{k-2}-\xi _{k-1}}{\xi _{k-1}-\xi _k}(x-\xi _k)+\xi _{k-1} & \xi
_k\leq x<\xi _{k-1} \\ 
\frac{x-a}{1-a} & a\leq x\leq 1
\end{array}
\right. $ 
\noindent
with $\xi _k=\frac a{(k+1)^{\frac 1{z-1}}}$ ,$\ k\in N,z \in {\bf R}, \ z\geq 2$.
% this is thepiecewise linear version of the manneville map $T(x)=x+x^z\ (mod\ 1)$.
This is a P.L. version of the Manneville map $T(x)=x+x^z \ (mod\  1)$ (see fig 1), this example will be discussed more deeply in Section 4.
 In this example any neighborhood of the
origin $B_\epsilon =$ $[0,\epsilon )$ is subdivided in a sequence of
intervals $A_k=(\xi _k,\xi _{k-1}]$ and if $k>1$ then  $T(A_k)=A_{k-1}$. Let us
choose $x=0,$ then  $B(n,0,\xi _k)=[0,\xi _{k+n}).$ By this we find $%
r(0)=[J\left( \frac a{(n+k+1)^{\frac 1{z-1}}}\right) ]$ that is a class of
infinitesimals corresponding to  power law decreasing sequences whit
exponent $\frac{-1}{z-1}$.  In other words the map $T(x)$  has power law
sensitivity to initial condition at the origin. 
In Section 4 we will see that while the sensitivity to initial condition at the origin is a power law,  for almost all other  points in $[0,1]$ we
have a stretched exponential sensitivity.
This example is important
in the applications and  will be  studied more deeply in section 4.

\subsection{orbit complexity}

Now we define our indicator of orbit complexity. In the philosophy of the algorithmic information content 
we define the complexity  of the orbit of $x$ as the asymptotic behavior of the quantity of information that is necessary to reconstruct the orbit, i.e. the asymptotic behavior (with respect to the variable $n$) of the length of the smallest program that can approximate $n$ steps of the orbit with its output (at accuracy $\epsilon$). As before we consider the behavior when $n$ goes to $\infty$ and the accuracy parameter goes to 0 .

To interpret the output of a calculation which is a finite string as a finite sequence in $X$ let us consider an interpretation function $I$ and 
 a total recursive surjective function.  
$${\cal Q}:\Sigma \rightarrow \Sigma ^* $$
where $\Sigma ^*$ is the set of finite sequences in $\Sigma $. Now let us
consider an  universal Turing machine ${\cal U}$,  for each program $p$  we  define $U(p)\in X^* $ (the set of finite sequences in $X$) as 

$$U(p)={I}( {\cal Q}({\cal U}(p)))$$

\noindent where $I$ is extended obviously to a map from the space $\Sigma ^*$ to $X^*$. $U_i(p)\in X$ is defined as the $i-$th point of $U(p)$ .
With this definition we can interpret the output of a calculation
as a finite sequence in $X$.
We  remark that given $\cal Q $ and a sequence of strings $s_1,...,s_n$
it is possible by an algorithm to    find a single string $s$ such that ${\cal Q} (s)=(s_1,...,s_n)$.

\begin{definition}
We define the algorithmic information content of the sequence $x,T(x),...,T^{n}(x)\in X^*$ at accuracy $\epsilon $ and with respect to the interpretation $I$
as:

$${\cal E}^I(x,n,\epsilon) =\min \left\{ |p|\ s.t. \ U(p)\in X^{n+1},
\stackunder {0\leq i\leq n} {max} (d(U_{i}(p),T^i(x)))<\epsilon\right\}. $$  
\end{definition}
\noindent As before we consider the behavior for $n\rightarrow \infty$ and define ${\cal E}^I(x,\epsilon) :X\times {\bf R}\rightarrow \frac {{\bf R}^*}{\simeq}$ as:
$${\cal E}^I(x,\epsilon)=[J({\cal E}^I(x,n,\epsilon))].$$

\begin{remark}$ {\cal E}^I(x,\epsilon)$ is a non increasing function
with respect $\epsilon$.
\end{remark}

Finally, like in the definitions of initial data sensitivity we consider
the behavior when $\epsilon$ goes to $0$ and   we define $ {\cal E}^I(x) :X\rightarrow {\cal R}$ as

\begin{definition}
The orbit complexity of $x$ with respect to the interpretation $I$ is defined as:
$${\cal E}^I(x)=\stackunder {\epsilon\in {\bf R}^+}{sup}   {\cal E}^I(x,\epsilon) .  $$
\end{definition}

We now give some example of different behaviors of ${\cal E}^I(x)$.
If $x$ is a periodic point it is easy to see that  ${\cal E}^I(x)\leq [J(log(n))]$.
By the results of \cite{brud} and \cite{io2} it follows (see also Section 1) that if a system is compact, ergodic and has positive Kolmogorov entropy then for almost all points we have ${\cal E}^I(x)=[J(n)]$.
We also remark that (when the space is compact)  this is the maximum over all the possible behaviors.
Indeed if $X$ is compact, for each $\epsilon$ there is a finite cover made
of balls with ideal center and radius $\epsilon$, then a program that follows $n$ steps of
the orbit of any point with the accuracy $\epsilon$ can be simply made by  listing $n$
centers of the cover, then, if $X$ is compact ${\cal E}^I(x)\leq [J(n)]$.
In section 4 we will study the complexity of the orbits of another, less 
trivial example.

%As before we remark that from lemma @ ${\cal E}^I(x)=\stackunder {\epsilon \rightarrow \infty}{lim}{\cal E}^I(x,\epsilon)$.

\begin{lemma}\label{lemma29}
If $I ,J$ are computable interpretation functions from the same computable structure: $I ,J\in {\cal I}$ then ${\cal E}^I(x)={\cal E}^J(x)$.
So the orbit complexity does not depend on the choice of the interpretation $I$ in the computable structure $ {\cal I}$ and we can define ${\cal E}^{\cal I}(x)={\cal E}^I(x)$ for some $I \in {\cal I}$.

\end{lemma}
{\em Proof.} 
Let us consider equivalent interpretations $I_1,I_2 $ and

\noindent  $U=I_1({\cal Q}({\cal U(*)}))$ as in the definition of orbit complexity. Let us suppose that we have a minimal length  program $p_k$ for the interpretation $I_1$  such that $\forall i<(k)$ we have  $d(U_{i}(p_k),T^i(x))<2^{-\lambda - 1}$, then there is a program
$p_k'$ for $I_2$  approximating the orbit of $x$ with accuracy $2^{-\lambda}$ and $|p_k'|<|p_k|+c$. The program $p_k'$  runs $p_k$  finding strings $s_i$ such that
$I_1(s_i)=U_i(p_k)$,  then using the equivalence between $I_1$ and $I_2$ it finds strings $z_i$
such that $d(I_2(z_i),I_1(s_i))<2^{-\lambda -1}$ (Remark \ref{remark10}) by these strings it is easy to see how $p_k'$ can approximate the orbit of $x$ with accuracy $2^{-\lambda}$. It follows
that ${\cal E}^{I_2}(x,2^{-\lambda})<{\cal E}^{I_1}(x,2^{-\lambda +1}) $, then we have ${\cal E}^{I_2}(x)\leq{\cal E}^{I_1}(x) $ and exchanging $I_1$ with $I_2$
we obtain the opposite inequality. $\Box$

 If \( X \) is compact then the orbit complexity does not depend not even on the
computable structure.

\begin{theorem}\label{th30} If $X$ is compact,
if \( I \) is a computable interprepretation and \( J \) is another interpretation
function (not necessarily computable) then 
${\cal E}^I(x)\geq{\cal E}^J(x)$.
\end{theorem}
{\em Proof.} 
Let $\epsilon > 0$, $s_1,...,s_k \in \Sigma $   be a finite set of strings such that $B_{\epsilon}(I(s_1) ),...,B_{\epsilon}(I(s_k) )$ is a cover of $X$. The set of strings is finite because $X$ is compact. 
It is easy to see that there is an algorithm $A:{\bf N}\times {\Sigma } \rightarrow \{1,...,k\}$ such that $A(i,s)=m$ implies that the $i-th$ point of  $I({\cal Q}(s))\in B_{\epsilon}(I(s_m))$.
That is: the algorithm gets a string and a natural number and outputs 
a set $B_{\epsilon}(I(s_m))$ of the cover in which the $i$-th point of the interpretation of the string as a sequence in $X$ is contained.
The algorithm calculates the distance between the $i-th$ point of  $I({\cal Q}(s))\in B_{\epsilon}(I(s_m))$ and $I(s_z)$ for all
$z\in\{1,...,k\} $ with accuracy $\frac \epsilon 2$, until it finds an $s_m$ such that $d(U_i(s),I(s_m))<\frac \epsilon 2$ this is possible because $I$ is
a computable interpretation.

Now let us consider the interpretation $J$. Even for the interpretation $J$ there is a finite set $\{ s'_1,...s'_{k'}\}$ such that $B_{\epsilon}(J(s'_1)) ,...,B_{\epsilon}( J(s'_{k'}))$ is a cover of $X$. Now let us consider a function
$G:\{1,...,k\}\rightarrow \{1,...,k'\}$ such that $G(i)=j$ if $I(s_i)\in B_{\epsilon}(J(s'_{j}))$. Being a function between finite sets $G$ is a recursive function.

Now let $p$ be a minimal length program that allows to follow the orbit 
of $x$ for $n$ steps with accuracy $\epsilon$ and interpretation $I$, that is 
$$\stackunder {0\leq i\leq n} {max} (d(U_{i}(p),T^i(x)))<\epsilon, |p|={\cal E}^I(x,n,\epsilon).$$
For each $i$  by calculating  $A(i,{\cal U}(p))$ we can find  an $m$ such that   
$T^i(x)\in B_{2\epsilon}(I(s_m))$ and then by  function $G$ we can find a $j$ such that
$T^i(x)\in B_{3\epsilon}(J(s_j))$.

Summarizing, this procedure allows (given the program $p$) to calculate
a sequence of strings $s_{j_1},...,s_{j_n}$ such that we can follow the orbit of $x$ with the 
interpretation $J$, for $n$ steps and accuracy $3\epsilon $.
This implies that $ |p|+c \geq {\cal E}^J(x,n,3\epsilon)$ where $c$ is the length of the above procedure and does not depend on $n$. From this we have ${\cal E}^I(x,\epsilon)\geq {\cal E}^J(x,3\epsilon)$ and ${\cal E}^J(x)\leq {\cal E}^I(x)$.
$\Box $

From the above theorem we see the curious fact that  in the compact case the orbit 
complexity reaches its maximum over all interpretations at a computable interptetation (a sort of Kolmogorov-Sinai theorem if we keep in mind the parallelism between orbit complexity and entropy) and the orbit complexity with respect to a computable structure does not depend on the choice of the computable structure ( if some computable structure exists on the space).
Moreover, all this is true independently of the properties of $T$.

\begin{corollary}
If $X$ is compact, if $I$ and $J$ are computable interpretations (not necessarily from the same
computable structure) then ${\cal E}^I(x)={\cal E}^J(x)$.
\end{corollary}

The orbit complexity is invariant for constructive isomorphisms of dynamical systems over non compact spaces, it  stated in the following propositions.
As before we remark that  if the space is  compact constructivity is not required. 
We omit the proofs that are similar to the previous ones.
\begin{theorem}
If $(X,d,T)$, $(Y,d',T')$ are topological dynamical systems over metric spaces
with  computable structures  ${\cal I},{\cal J}$  and 
$f$ is  onto and it is a  morphism between $(X,d,{\cal I})$ and $(Y,d',{\cal J})$
such that the following diagram 
\begin{equation}\label{equisopra}
\begin{array}{rcccl}
\ & \ & f & \ & \ \\
\ & X & \rightarrow & Y \ \\
T & \downarrow & \ & \downarrow & T'\\
\ & X & \rightarrow & Y \ \\
\ & \ & f & \ & \ \\
\end{array}
\end{equation}
\noindent commutes,if $x\in X$ and $y=f(x)\in Y$ then 
 ${\cal E}^I(x)\geq {\cal E}^J(y). $

\end{theorem}

\begin{theorem}
If $(X,d,T)$, $(Y,d',T')$ are topological dynamical systems over compact metric spaces
with  computable structures  ${\cal I},{\cal J}$, if $f$ is an homeomorphism $X\rightarrow Y$
such that the  diagram \ref{equisopra}
commutes and $x\in X$ and $y=f(x)\in Y$ then 
 ${\cal E}^I(x)= {\cal E}^J(y). $

\end{theorem}

\subsection{Complexity of points}
Now we define a function ${ S}^I(x,\epsilon ):X\times {\bf R}\rightarrow {\bf R}$, the function is a measure of the complexity of the points of $X$. The function is non increasing and  measures how much information is necessary to approximate a given point of $X$ with accuracy $\epsilon$. Thus it is a function that does not depend on the dynamics.
  In \cite{io1} a definition of local entropy for points of metric spaces was based on this idea and connections between $S$ and the concept of dimension are shown. In particular $S(x,\epsilon)$ is related to the local dimension of $X$ at $x$. 
\begin{definition}If $I$ is an interpretation function, ${\cal U}$ an universal computer we define the information contained in the point $x$ with respect to
the accuracy $\epsilon$ as: 
\begin{equation} \label{1} 
{ S}^I(x,\epsilon )={\min} \left\{ |p|\ \ s.t.\ \ d(I({\cal U}(p)),x)<\epsilon\right\}.
\end{equation}
\end{definition}
The function $S$ depends on the interpretation $I$. In the following we will avoid to mention explicitly the superscript $I$ when it is clear from the context.
The function $S$ depends also on the choice of ${\cal U}$.
As stated in section 2.3 this function can be extended to a function $S^*:X\times {\bf R}^*\rightarrow {\bf R}^*$.
Unfortunately $S^*$ may be not compatible with the relation $\simeq$, 
for this reason we define 
 $\overline {S'}:X\times {\bf R}^*/\simeq \rightarrow {\cal R} $ as follows:
$$\overline S'(x,\alpha)=\stackunder  {a\in\alpha}{sup}[S^*(x,a)].$$ 
    \noindent   If $\alpha\in {\cal R} $ is an equivalence class then $\overline S'(\alpha)$  is the equivalence class of the supremum value of  $S^*(a)$  where $a$ ranges in the class $\alpha$. Since an equivalence class in $ R^*/\simeq$
does not change by the adding of a constant, then the function ${S'}$ does not depends more
on the choice of the universal computer $\cal U$ in the definition of $ S$.
In the same way we define $\underline {S'}:X\times {\bf R}^*/\simeq \rightarrow {\cal R} $ as:
$$\underline {S'}(x,\alpha)=\stackunder {a\in\alpha}{inf}[S^*(x,a)].$$
 
Finally we extend $\overline {S'}$ and $\underline {S'}$ to functions 
 $\overline {S},\underline {S}:X\times {\cal R}\rightarrow {\cal R} $ as follows:
if $a\in {\cal R}$ and $a=p(a_i)$ (the $a_i$ are in 
$R^*/\simeq$ ) then we define $$\overline S(x,a)=\overline p(\overline {S'}(x,a_i)), \ \underline S(x,a)=\overline p(\underline {S'}(x,a_i)) $$
\noindent this is well defined because $\cal R$ is closed by countable $inf $ and $sup$ and it  does not depend on the choice of $a_i$ in the class $a$.
Because $S$ is monotonic and then $(a_i)\approx (a'_i) $ implies $S'(x,a_i)\approx S'(x,a'_i) $.

If $I$ and $J$ are in the same computable structure $\cal I $ then the functions $\overline S^I$  and $\overline S^J$  are equal.   $\overline S^I$ does not depend on the choice of the interpretation in the computable structure.

\begin{lemma}
$\overline S$ and $\underline S$  are independent of the choice of  $I\in {\cal I}$, in other words $I,J\in {\cal I}$ implies $\overline S^I (x,\epsilon)=\overline S^J (x,\epsilon) $  and $\underline S^I (x,\epsilon)=\underline S^J (x,\epsilon) $.
\end{lemma}
{\em Proof.} The proof is very similar to the proof of Lemma \ref{lemma29} and we omit it. $\Box $
% \begin{remark}
%???If X is finite dimensional , I is computable then $\overline S=\underline S$.
%ricoprimento con palle???
%If $X$ has finite $\epsilon$ capacity then $ \overline S=\underline S$ 
%\end{remark}
%@SUL BLOCCOtroviamo un ricoprimento di raggio $\epsilon$ ma tale che i centri siano una 
%$\epsilon / 2 $ rete allora il numero dei centri si puo stimare e quindi 
% un punto sara' dato dal numero del suo centro.non serve la compattezza?@ 
%!!!@@@!!!
%$N_\delta (F)=$ the smallest number of balls of radius $\delta $ that cover $
%F.$
%
%$N_\delta ^{\prime }(F)=$ the largest number of disjoints balls of radius $
%\delta $ with centers in $F.$
%
%\underline{$\dim $}$_B(F)=\stackunder{\delta \rightarrow 0}{\lim \inf }\frac{
%\log N_\delta (F)}{-\log \delta }=\stackunder{\delta \rightarrow 0}{\lim
%\inf }\frac{\log N_\delta ^{\prime }(F)}{-\log \delta }$
%
%$\dim _H(F)\leq $\underline{$\dim $}$_B(F)$

\noindent We remark that if $X= {\bf R}^n$ or $X$ is a finite dimensional manifold
then $\overline S=\underline S$.

\begin{remark} If $(X,d,{\cal I)}$ is a metric space with computable structure ${\cal I}$ and the lower box counting dimension \footnote{See e.g. \cite{falconer}.} of $X$ is finite : \underline{$\dim $}$_B(X)=d$ then
$S(x,\delta )\leq -d\log \delta +C$ where $C$ is a constant not depending on 
$x$ and $\delta $. Hence $\overline{S}$ and \underline{$S$} coincides for
all the points of $X.$
\end{remark}

\noindent The proof follows from the observation that if $X$ is  finite
dimensional, then the minimum number $n_\epsilon$ of balls in a cover of $X$  with radius $\epsilon$ is such that $n_\epsilon \sim \epsilon ^{-d}$. By the computable structure
we can construct the centers of a suitable  cover with  $n_\epsilon \sim \epsilon ^{-d}$   and obtain that each point of $X$ is approximated with accuracy $\epsilon$ by indicating a particular center of the cover, 
which costs $\leq log(n_{\epsilon})+c\leq -d\log(\epsilon)+C$  bits, where $c$ represents the length of the procedure that construct the centers of the suitable $\epsilon-$cover.

%@!{\em Proof.} We describe a procedure $A(m,\delta )=\{s_1,...,s_{k(m)}\}$
%that finds the strings indicating the centers of a collection of disjoint
%balls $B_\delta (_{I(s_1)}),...,B_\delta (_{I(s_{k(m)})})$ with centers $
%I(s_1),...,I(s_{k(m)})$ and radius $\delta $ such that $\stackunder{m\in N}{
%\cup }B_{2\delta }(_{I(s_1)}),...,B_{2\delta }(_{I(s_{k(m)})})$ is a cover
%of $X$ .
%
%The procedure runs inductively as follows: let us suppose to have found $n$
%centers $x_1=I(s_1),...,x_n=I(s_n)$ . To find the $n+1$ -th center $
%x_{n_{+1}}=$  $I(s_{n+1})$ the procedure will test all the strings in $
%\Sigma $ until it finds a string $s_{n+1}$ such that $\forall i\leq
%n\,d(x_i,x_{n_{\delta +1}})\geq \delta $ .
%
%The procedure stops at the $m$-th center. If $m\rightarrow \infty $ the
%procedure finds a set of $n_\delta $ centers with the above prescriptions.
%
%Since $n_\delta \leq N_\delta ^{\prime }$ then for each $\varepsilon $
%eventually with respect to $\delta $ we have $n_\delta \leq \delta
%^{-d-\varepsilon }.$
%
%Then each point of $X$ can be indicated with accuracy $2\delta $ by
%calculating the procedure $A(m_x,\delta )$, where $m_x\leq $ $n_\delta \leq
%\delta ^{-d-\varepsilon }$
%
%then $S(x,2\delta )\leq C+\log (m_x)\leq C+\log (\delta ^{-d-\varepsilon
%})\leq C-C^{\prime }\log (\delta ).\Box $
%
%@Bisogna dire prima che se $S(x,\delta )\leq -d\log \delta +C$ allora  $
%\overline{S}$ and \underline{$S$} coincides @

\begin{lemma}\label{lemma37}
If $x\in {\bf  R}^n$, $\mu $ is the Lesbegue measure on ${\bf R}^n$, if $\epsilon \rightarrow 0$ then for $\mu$-almost all
$x\in {\bf R}^n$, $S^I(x,\epsilon)=-nlog (\epsilon)+o(log(\epsilon))$ and if 
$a=\stackunder \epsilon {inf} [J(a_{n,\epsilon})]$ then  $\underline S(x,a)=\overline S(x,a)=\stackunder \epsilon {inf}[J(log(a_{n,\epsilon}))] $. 
\end {lemma}
{\em Proof}. We prove that for almost all $x\in {\bf R}^n$ $\stackunder {\epsilon \rightarrow 0 }{lim}\frac{S^{\cal I}(x,\epsilon )}{-nlog(\epsilon )}= 1$. Theorem 12 of \cite{io1} states that the set 
\begin{equation}\label{mistero}
W^d=\{ x\in {\bf R}^n\ s.t. \stackunder {i\rightarrow \infty }{liminf}\frac{S^{\cal I}(x,2^{-i})}{i}\leq d\}
\end{equation}
 has Hausdorff dimension less or equal than $d$.
This implies that  if $\epsilon =2^{-\gamma}$ then  $\stackunder {\epsilon \rightarrow 0 }{liminf}\frac{S^I(x,\epsilon )}{-nlog(\epsilon )}=
\stackunder {\gamma \rightarrow \infty }{liminf}\frac{S^I(x,2^{-\gamma} )}{n\gamma}  $. If $\stackunder {\epsilon \rightarrow 0 }{liminf}\frac{S^I(x,\epsilon )}{-nlog(\epsilon )}\leq 1$ then $x\in W^n$  because $S(x,2^{-\gamma-1}) \leq S(x,2^{int(-\gamma)})\leq S(x,2^{-\gamma+1})$\footnote{$int(r)$ denotes the integer part of $r$}. This implies that the set of the $x$ s.t.  $ \stackunder {\epsilon \rightarrow 0 }{liminf}\frac{S^I(x,\epsilon )}{-nlog(\epsilon )}<1$ is included in the set  $  \stackunder {d<n} \cup
W^d  $ that is a $0$ measure set.
 To prove the other inequality it is enough to remark that each $x\in {\bf R}$ can be approximated with accuracy $\epsilon$ by specifying its first $-int(log(\epsilon))$ digits, if $x\in {\bf R}^n$ we need $-nlog(\epsilon)$ digits to explicit 
the $n $ coordinates.$\Box$

\subsection{Initial data sensitivity and orbit complexity} 

Now we are ready to state the first proposition linking orbit complexity
to initial data sensitivity.

\begin{proposition}\label{quellagrossa}
 If $(X,T)$ is a dynamical system on a space with a computable structure $\cal I$, $I\in {\cal I}$ and $T$ is constructive. There are constants $c_1 $ and $ c_2$ such that For all $x\in X,n\in {\bf N},\epsilon \in {\bf R}^+$ 
\begin{equation}\label{quella1}
{\cal E}^I(x,n,2\epsilon)<{S}^I(x,r(x,n,\epsilon))+log (n)+c_1
\end{equation}

\begin{equation}\label{quella2}
{ S}^I(x,R(x,n,3\epsilon))\leq {\cal E}^I (n,x,\epsilon)+c_2.
\end{equation}
\end{proposition}

\noindent {\em Proof of \ref{quella1}.}
We will see that there is a program $p$ such that  
\noindent $d(U_{i}(p),T^i(x))<2%
\epsilon ,\forall i\leq n$ and $|p|\leq S(x,r(x,n,\epsilon)   )+\log (n)+c$.
If we have a program $p_0$ such that $d(I({\cal U}(p_0)),x)\leq r(x,k,\epsilon )$ then $ \forall i\leq k \
d(T^i(x),T^i(I({\cal U}(p_0)))<\epsilon.$

The idea is that by constructivity if we have the string  $s_0={\cal U}(p_0)$
we can follow the orbit of $I(s_0)$ by an algorithm $A(s_0,k,\epsilon )$
(see Lemma \ref{lemma12}). The program $p$ will codify the following procedure:

1) run the program $p_0$ and compute $s_0={\cal U}(p_0)$

2) compute $s_i=$ $A(s_0,i,\epsilon ),\forall 1\leq i\leq k$

3) compute the single string $s$ such that $Q(s)=(s_0,...,s_k).$

The length of this program will be a constant (the above stated procedure)
plus $\log (k)$ (the length of a binary representation of $k$) plus the
length of $p_0$. If $p_0$ was supposed to be the shortest program such that $%
d(I(U(p_0)),x)<r(x,k,\epsilon )$ then its length is the value of $S( x,r(x,n,\epsilon)  )$
and the first part of the statement is proved.

{\em Proof of \ref{quella2}.}
Let $e\in {\bf N}$ such that $2^{-e}<\epsilon$. If we have a program $p$  such that $d(T^i(x),U_{i}(p_0))<\epsilon $ for $0<i<k$
we can find a string $s$ s.t.  $I(s)\in B(x,n,3\epsilon )$ with the following procedure:

1) By $p_0$ compute the number $k$ and the strings $s_{0},...,s_{k}$ such that  $%
(s_{0,...,}s_k)=Q({\cal U}(p_0))$

2) for each $c\in \Sigma $ do the following things: \{
compute $A(c,i,2^{-e} )$ for each $0\leq i\leq k$ , if for all $0\leq i\leq k$ $D(A(c,i,2^{-e}
),s_i,e +2)<2^{-e} $ then $s=c$ and stop the procedure.
\}

The procedure must stop because of the density of the image of  $I$. At
some time the step 2) will be computed with a $c $ such that $I(c)\in B(x,n,2^{-e-1})$ and
this string will verify 2). On the other hand if we find a   $c$ that stops
the procedure then it is easy to see that $I(c)\in B(x,n,3\epsilon ).$ This
will implies that $d(x,I(c))<R(x,k,3\epsilon )$. Summarizing we have described a procedure that starting from a program $p_0$ outputs a string $s$ such that $d(I(s),x)<R(x,k,3\epsilon)$. The code for this procedure will be a program containing $p_0$ and its length will be $p_0+C$ where $C $ represents the length of the code for the above procedure which does not depend on $x$ and $k$, and the statement is proved.$\Box$

From the previous statement we obtain a relation between the indicators of
orbit complexity and sensitivity.

\begin{theorem}\label{quellogrosso}
If $(X,T)$ is a dynamical system on a space with a computable structure $\cal I$ and $T$ is constructive. For all $x\in X$

\begin{equation}\label{e4}
{\cal E}^{\cal I}(x)\leq max(\overline {S}^{\cal I}(x,r(x)),[J(log (n))])
\end{equation}
\begin{equation}\label{e5}
\underline{S}^{\cal I}(x,R(x))\leq {\cal E}^{\cal I} (x).
\end{equation}
\end{theorem}
{\em Proof.}
If we apply the homomorphism $J$ to equation \ref{quella1} we obtain $$J({\cal E}^I(x,n,2\epsilon))\leq J(S^I(x,r(x,n,\epsilon)))+J(log (n))$$, then 
$J({\cal E}^I(x,n,2\epsilon)) \leq {S^*}^I(x,J(r(x,n,\epsilon)))+J(log (n))$, and considering the equivalence classes:
$${\cal E}^I(x,2\epsilon)\leq max(\overline{S'}^I(x,r_{\epsilon}(x)),[J(log (n))]).$$
This is true for each $\epsilon$, then 
${\cal E}^I(x)\leq max( \overline {S}^I(x,r(x))+[J(log (n))])$. As proved before all this equivalence classes does not depend on the choice of $I\in {\cal I}$ and we have
Equation \ref{e4}.
In the same way we can obtain equation \ref{e5}.$\Box$

By Lemma \ref{lemma37} for almost all points in ${\bf R}^n$ the function $S$ is the logarithm, this, combined with proposition \ref{quellagrossa} implies the following formulas:

\begin{theorem}If $T:{\bf R}^n\rightarrow {\bf R}^n$ is constructive on ${\bf R}^n$ with the standard computable structure,  for almost all $x\in X$

\begin{equation}{\cal E}^{\cal I}(x)\leq max(\stackunder{\epsilon \in {\bf R}^+}{inf}([J(log(r(x,n,\epsilon))]),[J(log (n))])
\end{equation}
\begin{equation}\stackunder {\epsilon \in {\bf R}^+}{inf}[J(log(R(x,n,\epsilon )))]\leq {\cal E} (x)
\end{equation}
\end{theorem}

As a corollary of Theorem \ref{quellogrosso} we can obtain the following interesting result:
the set where the sensitivity to initial conditions is more than exponential
in all directions has $0$  Hausdorff dimension.
\begin{theorem}\label{th41}
If $(X,T)$ is a dynamical system on a compact metric space with a computable structure $\cal I$ and $T$ is constructive. Then the set $$\overline {exp}=\{ x\in X \ s.t. \forall h\in {\bf R}^+\ R(x)<[J(2^{-h n})]\}$$ has zero Hausdorff dimension.
\end{theorem}

{\em  Proof. } %Let $a_i:{\bf N}\rightarrow {\bf R} $ be the sequence given by $a_i=i$.
 Let us consider a point $x$ such that $R(x)<[J(2^{-h n})] $ by theorem \ref{quellogrosso} we know that $\underline{S}^{\cal I}(x,R(x))\leq {\cal E}^{\cal I} (x)$
 since $X$ is compact we have ${\cal E}^{\cal I} (x)\leq [J(n)]$.
Then $\underline{S}^{\cal I}(x,R(x))\leq [J(n)]$. Since $\underline{S} $
is a non increasing function then $\underline{S}^{\cal I}(x,[J(2^{-h n})])\leq \underline{S}^{\cal I}(x,R(x))\leq [J(n)]$.
Then by definition $\underline{S}^{\cal I}(x,[J(2^{-h n})])=
 \stackunder {\alpha \in [J(2^{-h n})]} {inf}[S^*(x,\alpha)]\geq[S^*(x,J(2^{-(h-\epsilon) n}))]$ for some small $\epsilon$ and then $[S^*(x,J(2^{-(h-\epsilon) n}))] \leq [J(n)]$, by this, setting $h'=h-\epsilon$ it follows that there is a bounded 
constant $c_1 \neq 0$  such that $ S^*(x,J(2^{-h' n}))\leq c_1  J(n)$.

Let us consider the following set $A^{d}=\{x\in X|J(\frac{S(x,2^{-i})}i)\leq
d \}$ , since for each sequence $b_i$ we have $J(b_i)\geq \lim \inf (b_i)$
then $A^d\subset W^d$ where $W^d$ is the set defined in eq. \ref{mistero}. Since the
Hausdorff dimension of $W^d$ is greater or equal than $d$ then also $\dim
_H(A^{d})\leq d$ . Now let us consider the set $$A^{d,h'}=\{x\in X|J(\frac{
S(x,2^{-h'i})}i)\leq d\}$$.

\noindent If $x\in A^{d,h'}$ , let us set   $k=h'i$ and let us consider $\frac{
S(x,2^{-k})}{\frac k{h'}}$, if $J(\frac{S(x,2^{-h'i})}i)\leq d$ then $
\stackunder{k\rightarrow \infty }{\lim \inf }(\frac{S(x,2^{-k})}{\frac k{h'}}
)\leq J(\frac{S(x,2^{-h'i})}i) \leq d$ because $\frac{S(x,2^{-h'i})}i$ is a subsequence of $\frac{
S(x,2^{-k})}{\frac k{h'}}$ (Axiom \ref{ax16}). Then $x\in W^{\frac d{h'}}.$ This implies that if $x\in \overline{exp}$ then $x\in \stackunder {h'} {\cap}W^{\frac d {h'}}$ which has $0$ Hausdorff dimension (again by \cite{io1} Theorem 12).
$\Box
$

\section{Applications to the Mannevile maps.}

In this section, in order to give a non trivial example of application of the theory exposed in the previous sections we present some example of weakly chaotic dynamics.
We construct a class of examples of dynamical systems over the unit interval with stretched
exponential sensitivity to initial conditions and information
content of the orbits that increases as a power law.
We precise that the maps $T_z$ we are going to study are not weakly chaotic
in the sense of \cite{smital} (zero topological entropy), conversely they have positive topological entropy. In this examples however for almost all the points (for the Lesbegue measure) the dynamics are  weakly chaotic (low orbit complexity, low initial data sensitivity). Then we can say that the system 
is weakly chaotic with respect to the Lesbegue measure.

The examples are {\em piecewise linear}  version of the Manneville map $T:[0,1]\rightarrow [0,1]$ defined as 
$T_z(x)=x +x^z \ (mod \ 1), \ z \in (1,\infty) $.
The so called Manneville map comes from the theory of turbulence.
 It was introduced  in \cite{Manneville} as an extremely
simplified model of intermittent behavior in fluid dynamics,
then the map was studied and applied in other areas of the physics (for example \cite{Grigolini},\cite{toth},\cite{pollicot}).

The first study of the mathematical features of the Manneville map was done by Gaspard and Wang in \cite{Gaspard}. However in our opinion In their paper some steps of the proofs were
 difficult to understand and some others were not rigorously formalized. 
In the following we outline the construction done in \cite{Gaspard} for the 
study of the complexity of the  piecewise linear  Manneville maps by the theory of recurrent events \cite{feller}. Then we  prove the 
main features of  this important class of dynamical systems by the theory
exposed in the previous sections. Another  study of the Manneville
map was done by C. Bonanno in \cite{bonanno} where the dynamics was studied also from
a topological point of view. Notations: if $(a_i),(b_j):{\bf N}\rightarrow {\bf R}$ are real sequences, in the following we will write $(a_i) \underline {\simeq} (b_j)$ if and only if $\frac {a_i}{b_i}$ is bounded, we also write $(a_i)  {\sim} (b_j)$ if and only if $\stackunder i\lim \frac {a_i}{b_i}=1.$

Let  \(\epsilon=( \epsilon _{k}):N\rightarrow [0,1] \) be a monotone sequence, such that
\(\stackunder {k\rightarrow \infty }{lim}\epsilon _{k} \)\( =0. \) Let \( T_\epsilon:[0,1] \)\( \rightarrow [0,1] \)
be defined by:

\( T_\epsilon (x)=\left\{ \begin{array}{cc}
\frac{\epsilon _{k-2}-\epsilon _{k-1}}{\epsilon _{k-1}-\epsilon _{k}}(x-\epsilon _{k})+\epsilon _{k-1} & k>0,\epsilon _{k}<x\leq \epsilon _{k-1}\\
\frac{x-\epsilon _{0}}{1-\epsilon _{0}} & \epsilon _{0}<x\leq 1
\end{array}\right.  \)

\noindent To each sequence \( \epsilon _{i} \) is then associated a piecewise linear map \( T_\epsilon \) (see fig. 1) and a dynamical system $([0,1],T_\epsilon)$.

In order to apply the theory of recurrent events, we associate to our dynamical system \( ([0,1],T_\epsilon) \)  a stochastic process 
\( X _{i} \).The process is defined on the probability space \( (\Omega ,\mu ), \) where 
\( \Omega =[0,1] \) and \( \mu  \)
is the Lesbegue measure as follows. Let us consider the sets $A_i, i\in {\bf N}$: \( A_{0}=(\epsilon _{0},1],...,\, A_{i}=(\epsilon _{i},\epsilon _{i-1}] \).
Let \( f:{\bf R}\rightarrow {\bf N}  \) be given by \( f(x)=n\Leftrightarrow x\in A_{n} \)
(\( f \) associates to each point \( x \) the index of the set \( A_{i} \)
in which \( x \) is included). The associated stochastic process is given by the random variables \( X_{i}:\Omega \rightarrow {\bf R} \) , given by
\( X_{i}(\omega )=f(T_{\epsilon} ^{i}(\omega )) \). 
As it was remarked in \cite{Gaspard} the process \( X_{i} \) is a Markov chain with transition matrix 

\( \left( \begin{array}{ccccc} 
p_{0} & p_{1} & p_{2} & p_{3} & \cdots   \\
1 & 0 & 0 & 0 & \cdots    \\
0 & 1 & 0 & 0 & \cdots    \\
0 & 0 & 1 & 0 & \cdots    \\
\vdots  & \vdots  & \vdots  & \vdots  & \ddots  \\
 
\end{array}\right)  \) where \( p_{i}=\mu (A_{i})=\mu (X_{k}=i|X_{k-1}=0). \) 
\begin{figure}
\begin{center}
\includegraphics[height=7cm]{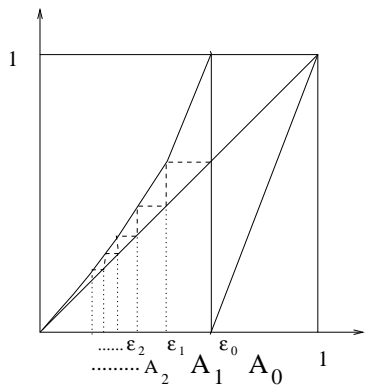}
\caption{\label{xsosx}The  map  $T_z$ and the partition $A_i$}
\end{center}
\end{figure}
Summarizing:     we constructed a 
 family of dynamical systems (one for each infinitesimal sequence). To each one of them it is associated a Markov chain. The statistic behavior of the Markov chain can be studied by the theory of recurrent events \cite{feller} and will give information on the dynamics. 
The family $([0,1],T_\epsilon)$  is a large family of dynamical systems with different chaotic behavior.
In the following we will study the systems in which \( \epsilon _{i} \sim \frac1{i^\alpha}\). Such dynamical systems give rise to stretched exponential
initial data sensitivity and power law orbit complexity.
However if \( \epsilon _{i} \) decreases slower than a power law (\( \epsilon _{i}\sim \frac {-1}{log(i)}    \) e.g.) we will have a variety of other possible
behaviors of orbit complexity and initial data sensitivity that will be not
studied here.

Now we consider a subset  of the above  family of  systems, which we consider as P.L. versions of the Manneville
map.
Let $z\in (2,\infty )$, let us consider the map $T_z$ associated by the above construction to  the sequence  $$\epsilon _k = 1 - \sum_{i\geq k}\frac{1}{(i+1)^{\frac{z}{z-1}}}. $$

\noindent The map $T_z$ is then  associated to a Markov chain with  transition probabilities  \( p_0=1-\epsilon _0, \ p_{i}= \frac{1}{(i+1)^{\frac{z}{z-1}}} \    \forall i>1  \). %%%%%%@Moreover we remark that if $z$ is a constructive number (a number that can be approximated at any accuracy with finite information) the map $T_z$, defined on the space $[0,1]$ with the standard computable structure  is constructive and then we can apply  Proposition \ref{quellagrossa} to $T_z$.

The theory of recurrent events  can now give us information on the dynamics of the Markov chain. Let us consider the random variable $N_n(x)$ that is, 
the number of times that the event $A_0$ occurs until the time $n$: 
$N_n(x)=\#\{i|X_i(x)=0\}=\#\{i|T^i(x)\in A_0\}$. Theorem 10 of \cite{feller} applied to our
Markov chain says that if $1-\sum_{i\geq k}p_i \sim Ak^{-\alpha}$
then $E(N_n)\sim C n^{\alpha}$. Since for the Markov chain associated to  $T_z$ we have $1-\stackunder{i\geq k}\sum p_i\sim \frac 1{z-1}k^{\frac {1}{z-1}}$ then $E(N_n)\sim C n^{ \frac {1}{z-1}}$.

Since the function $T_z$ is not continuous at $\epsilon _0$ our theory is not directly applicable to the dynamical system $([0,1],T_z)$. The following lemma and its proof shows
how to extend our theory to discontinuous dynamical systems. 
A real number is constructive if it can be approximated at any accuracy 
by an algorithm.
\begin{definition}A number $z\in {\bf R}$ is said to be constructive if there
is an algorithm $A_z(n):{\bf N}\rightarrow {\bf Q}$ such that $A_z(n)=q $ implies
$|q-z|<2^{-n}$. 
\end{definition}
The rational numbers, the algebraic numbers are constructive, and so are 
all the numbers that can be explicitly used for numerical purposes,
for example $\pi$ and $e$ are constructive.

\begin{lemma} If $z$ is constructive. For the dynamical system $([0,1],T_z)$,  for all $x\in [0,1]$ equation \ref{quella1} and \ref{quella2} holds.
\end{lemma}
{\em proof.}
Let us consider the set $X=[0,1] - \{x|\exists k\in {\bf N},T_z^k(x)=\epsilon _0\}$. Let us consider an  interpretation function $I'$ on $X$ as follows: if $I$ is defined as in equation \ref{I}
 then  $I':\Sigma \rightarrow X$ is  given by
$I'(s)=I(s)+\pi$. $I'$ is a computable interpretation because $\pi $ is constructive and its image is in $X$ because $\pi$ is transcendent. If $z$ is constructive then it is easy to see that $T_z$ is constructive on $(X,d,\cal {I})$
where $I'\in \cal{ I}$ then we can apply Proposition \ref{quellagrossa}  to $(X,T_z,{\cal I})$.
On the other hand, since the inclusion $i:X\rightarrow [0,1]$ is isometric the complexity of an orbit in $X$ is equal to the complexity of the corresponding orbit in $[0,1]$: $\forall x \in X$ $ {\cal E}^{I'}(([0,1],T_z),x,n,\epsilon)=
{\cal E}^{I'}((X,T_z),x,n,\epsilon)$. If conversely 
$x\in [0,1]-X$ then ${\cal E}^{I'}(([0,1],T_z),x)=[J(log(n))]$ because
the orbit of $x$ converges to a fixed point. Then eq. \ref{quella1} and \ref{quella2} holds for all points of $[0,1]$ for the interpretation $I'$. 
 Now, since $[0,1]$ is compact 
(see theorem \ref{th30}) we have that the orbit complexity does not depend on $I'$ and the statement is proved.
 $\Box$

%***
%
%\begin{lemma}
%If there is a set $W\subset [0,1]$ of positive measure such that $\forall x \in W N_n(x)= o(n^\alpha) $ then $[J(E(N_n))]\geq [J(n^\alpha)]$.
%

%\end{lemma}
%
%{\em proof}If $\mu (W)=\epsilon $ then $E(N_n)\geq \epsilon E(N_n|W)$...
%
%

%\begin{lemma}
%If $E(N_n)\sim n^\alpha $ then there is a set $W$ with positive measure
%such that $ \forall x \in W [J(N_n(x))]\geq [J(n^\alpha)] $ 
%\end{lemma}
%

%\begin{lemma} 
%1+2=$[J(E(N_n))]=sup\{ f\in {\cal R}|\exists W,\mu (W)>0, s.t. \forall 
%x \in W, [J(N_n(x))]=f\}$ 
%\end{lemma}

Now let us give an estimation of the initial data sensitivity
of the map $T_z$ by the behaviour of $N_n$.
Let $x,y\in [0,1], d(x,y)=\Delta x(0)=\epsilon_0$.
Since the derivative of $T_z$ exists for almost all points and it is greater than $1$, if $\epsilon _0$ is small enough (e.g. is such that $\epsilon _0\leq \frac {p_0}4$) then $d(T_z(x),T_z(y))\geq \epsilon _0$.
 If $\epsilon _0$ is small enough (e.g. $\epsilon _0$ is such that $\forall i \leq n , d(T^i_z(x),T^i_z(y))\leq \frac {p_0}4$) , then $d(T^n_z(x),T^n_z(y))=\Delta x(n)\geq \epsilon _02^{cN_n}$ where $c$ is a constant 
that depends on the derivative of $T_z$ in $A_0$ and $A_1$.

If $d(x,y)=\epsilon _0=\epsilon 2^{-c N_n}$ then $\Delta x(n)\geq \epsilon _0 2^{c N_n(x)}\geq \epsilon.$ This implies that $y\notin B(n,x,\epsilon)=\{ y \in X :d(T^i(y),T^i(x))\leq \epsilon \ \forall \ 0\leq i \leq n\}$ (as defined in section 3.1)  and  symmetrically $y'=2x-y\notin B(n,x,\epsilon)$, then
$\epsilon _0 =\epsilon 2^{-c N_n(x)}\geq R(n,\epsilon,x)$ (and then $R(x)\leq [J(2^{-c N_n(x)})]$).

By this we can give an estimation of the complexity of the orbits of the Manneville maps. In analogy with the results of \cite{Gaspard} we are ready to state the following proposition 
about the asymptotic behavior of the average  orbit complexity in the Manneville map.  

\begin{proposition}
If $z$ is constructive. The average asymptotic behavior of the complexity of the orbits of $T_z$ is 
for each $\epsilon$  $$E({\cal E}^I (x,n,\epsilon))\underline {\simeq}n^\alpha $$ 
\noindent where $\alpha={\frac {1}{z-1}}. $
\end{proposition}
{\em Proof.}
Since we  have $R(n,\epsilon,x)\leq \epsilon 2^{-cN_n(x)}$  by proposition \ref{quellagrossa}
and proposition \ref{lemma37} 
$-log(R(x,n,3\epsilon))+o(log(R(x,n,3\epsilon))) \leq {\cal E}^I (x,n,\epsilon)+c_1$ for almost all $x\in [0,1]$,
then ${\cal E}^I (x,n,\epsilon)\geq c N_n (x)+o( N_n (x))+const+log (\epsilon)$ and
$E({\cal E}^I (x,n,\epsilon))\geq c E(N_n (x))  +E(o( N_n (x)))  +const+log(\epsilon)$.
But since $E(N_n (x))\sim n^\alpha$ by the theory of recurrent events,
 we have that for all $\epsilon$, $n^\alpha=O(E({\cal E}^I (x,n,\epsilon)))$ (that implies $[J(E({\cal E}^I (x,n,x,\epsilon)))]\geq [J(n^\alpha)]$ with our notations).

Now let us give an estimation from above of the average orbit complexity
in the Manneville maps. As before we see that the number $N_n(x)$ will
give the main part of the complexity of the orbit of $x$.

Let us consider an interpretation $I$ for the standard computable
structure in $[0,1]$ and  a minimal cover $B_\epsilon(I(s_0)),...,B_\epsilon(I(s_l))$ of $[0,1]$ (i.e. a cover of $[0,1]$ such that each its proper subset 
does not cover $[0,1]$). The number $l$ of balls in this cover is then bounded
by $l\leq \frac 1 \epsilon$. Now let us consider the sets $A_i$ defined above.
We remark that if $ 1 - \sum_{i\geq k}\frac{1}{(i+1)^{\frac{z}{z-1}}}<\epsilon $ 
then $A_k\subset B_\epsilon(I(s_0))$. Now we remark that the symbolic dynamics 
of the point $x$ with respect to $\{A_i\}_{i \in {\bf N}}$ i.e. the sequence
of $A_{i_0},...,A_{i_n}$ such that $T_z^n(x)\in A_{i_n}$ it is determined by 
the recurrence times of the set $A_0$. The sequence $(A_{i_n})=A_{i_0},...,A_{i_k}$ must be such that if $i_n>0$ then $A_{i_{n+1}}=A_{{i_n}-1}$ and if
$i_n=0$ then $A_{i_{n+1}}$ can be any one of the $A_i$ (see fig 1). For example a possible sequence is $A_0,A_1,A_0,A_3,A_2,A_1,A_0,A_0,A_2,...$ such a string is determined by
the sequence $P$   of numbers representing the recurrence times of $A_0$, i.e. for the above example  $P=(0,1,3,0,2,...)$.  We remark that the string $P_n(x)=(t_1,...,t_{N_n(x)})$ representing the recurrence times of $A_0$ for $n$ steps of the orbit of $x$ contains $N_n(x)$ numbers.
Then the binary length of $P_n(x)$ is about $\sum^{N_n}_{i=1} log(t_i)$.
  
 An algorithm $A_\epsilon (n,m)$ to follow the orbit of $x$ with accuracy $\epsilon$ (i.e. such that if $n\leq m$ and $A_\epsilon (n,m)=b$
then $T_z^n(x)\in B_\epsilon (I(s_b))$ ) can be constructed in the  following way.
The algorithm contains a string $P_m(x)$ %for the simbolic  orbit of $x$ with respect to ${A_i}_{i \in {\bf N}}$. 
of the recurrence times with respect to $A_0$, with this string we can reconstruct the symbolic orbit of $x$ with respect to $\{A_i\}_{i\in {\bf N}}$ for $m$ steps, as described above.
 Moreover the program contains another string $Q_m(x)$
containing at most $(l-1) N_n(x)$ numbers, each one is less or equal than $l$.
The meaning of $Q_m(x)$ will be explained below.

 The algorithm starts with a pointer to the first number of $Q_m(x)$. By $P_m$ it calculates the set $A_{i_n}$ such that 
 $T_z^n(x)\in A_{i_n} $, if $i_n$ is such that $ 1 - \sum_{j\geq i_n}\frac{1}{(j+1)^{\frac{z}{z-1}}}\leq \epsilon$ then $A_\epsilon (n,m)=0$ because then $T_z^n(x)\in B_\epsilon (s_0)$. Else the algorithm outputs the number of $Q_m(x)$ indicated
by the pointer and then set the pointer on the next number.
In other words, when  $T_z^n(x)$ may not be in $ B_\epsilon (s_0)$ then 
the algorithm gets the number $q$ such that $T_z^n(x)\in B_\epsilon (s_q)$
from the list  $Q_m(x)$.  This can not be too expensive because the point
can be out of $ B_\epsilon (s_0)$ at most  $(l-1) N_n(x)$ times.
The total length of the program implementing  $A_\epsilon (n,m)$  for $m$ steps
of the orbit is then
less or equal than $ \sum^{N_n}_{i=1} log(t_i)+(l-1) N_n(x)log(l)+C$ where $C$ is constant
with respect to $m$, $l$ depends on $\epsilon $ but not on $m$.
 
The term $\sum^{N_n}_{i=1} log(t_i)$ can be estimated as follows: we remark that while
the random variables $t_i$ are independent and identically distributed they have no finite expectation 
when $p_i\sim \frac 1 {i^\alpha}$, $1<\alpha <2$ ($E(t_i)=\sum ip_i$) instead the random variable 
$log(t_i)$ has finite expectation, let us say $E(log(t_i))=\overline {t}$. Then $\sum^{N_n}_{i=1} log(t_i)=N_n(\frac {\sum^{N_n}_{i=1} log(t_i)}{N_n})$, by the law of large numbers we have that for each
$\delta>1$, for almost each $x$, eventually with respect to $i$ we have $\frac {\sum^{N_n}_{i=1} log(t_i)}{N_n}<\delta\overline {t}$.

This implies that for almost each $x$, eventually with respect to $n$  ${\cal E}^I (x,n,\epsilon)\leq c N_n (x)\delta\overline t +const $ and then
$E({\cal E}^I (x,n,x,\epsilon))\geq C E(N_n (x))+const$.
But since $E(N_n (x))\sim n^\alpha$
 we have that for all $\epsilon$, $E({\cal E}^I (x,n,\epsilon))=O(n^\alpha)$. $\Box$

The above Theorem is an estimation of the average orbit complexity of the map $T_z$,
to  show how our theory can be applied to prove rigorously 
the statements of \cite{Gaspard},
however stronger results can be proved. By \cite{Gaspard} page 4592 eq. 2.9 (which follows from \cite{feller} theorem 7, page 106) we have that there exists a constant $A$ such that $$ \mu(\{N_n\geq \frac {n^\alpha}{Ax^\alpha}\})\rightarrow G_\alpha (x)$$ \noindent where $\alpha=\frac 1{z-1}$ and $G_\alpha$ is the Levi stable distribution law with parameter $\alpha$.
It follows that $N_n(x) \underline \simeq C n^{\alpha}$ for almost all points in the interval. 
 %!!!questa da una diseq. l'altra
%per assurdo se ci fosse un insieme di misura maggiore di zero con
%potenza maggiore farebbe media da solo
 From this, by the same proof as above it follows the   pointwise estimation: 

\begin{theorem}With the same notations as above,
For almost all $x\in [0,1]$, ${\cal E} (x,n,\epsilon))\underline \simeq n^\alpha $. 
\end{theorem}

Similar results are obtained by \cite{bonanno} using different techniques.

\section{Numerical experiments}
We want to remark that while the information content of an orbit 
(as it is defined in this work)
is not computable (the algorithmic information content of a string
is not a computable function) there is the possibility 
to have an empirical estimation of the quantity of information  by the use of
data compression algorithms. If instead to measure the information
contained in a string by its algorithmic information content
we consider as 'approximate' measure of the information content of the string the length of the string after it is compressed 
by a suitable coding procedure we obtain a computable notion of orbit complexity. In the positive entropy case the computable orbit complexity is a.e. equivalent to the previous
one \cite{io3}. 
Such a definition of computable orbit complexity allows numerical
investigations about the complexity of unknown systems. Unknown systems underlying  for example some given time series or experimental datas.

The existence of a computable version of the orbit complexity
motivates from the applicative point of view the study of the  orbit complexity itself and its relations between the other  measures of  the chaotic behavior of a system.

In \cite{ZIPPO} and \cite{menconi} such numerical investigations are performed by directly measuring the complexity of the orbits
of the Manneville map and of the logistic map at the chaos threshold. The
results  agree with theoretical predictions and some interesting 
conjecture arises.

%But $\underline{S}(x,{\cal I},R(x))\leq {\cal E} (x,{\cal I})$
%and for almost all $x\in [0,1]$ $log(R(x))\leq {\cal E} (x,{\cal I})$, this implies that for almost all $x$ in $[0,1]$ we have $[log(\epsilon 2^{-c N_n})]\leq {\cal E} (x,{\cal I})$. That is $[{\cal E} (x,{\cal I})]\geq [N_n]$, by lemma 39 o theo 37...

\end{document}